\documentclass{article}
\usepackage{amsmath,amsthm}
\usepackage{graphicx,epstopdf,epsfig,multirow,epic,bm}
\usepackage{lineno,hyperref}
\usepackage{subfigure}
\usepackage{amsfonts}

\normalsize \rm
\begin{document}

\begin{center}
{\Large  \textbf { Scale-free tree network with an ultra-large diameter }}\\[12pt]
{\large Fei Ma$^{a,}$\footnote{~The author's E-mail: mafei123987@163.com. },\quad  Ping Wang$^{b,c,d,}$\footnote{~The corresponding author's E-mail: pwang@pku.edu.cn.} }\\[6pt]
{\footnotesize $^{a}$ School of Electronics Engineering and Computer Science, Peking University, Beijing 100871, China\\
$^{b}$ National Engineering Research Center for Software Engineering, Peking University, Beijing, China\\
$^{c}$ School of Software and Microelectronics, Peking University, Beijing  102600, China\\
$^{d}$ Key Laboratory of High Confidence Software Technologies (PKU), Ministry of Education, Beijing, China}\\[12pt]
\end{center}

\begin{quote}
\textbf{Abstract:} Scale-free networks are prevalently observed in a great variety of complex systems, which triggers various researches relevant to networked models of such type. In this work, we propose a family of growth tree networks $\mathcal{T}_{t}$, which turn out to be scale-free, in an iterative manner. As opposed to most of published tree models with scale-free feature, our tree networks have the power-law exponent $\gamma=1+\ln5/\ln2$ that is obviously larger than $3$. At the same time, "small-world" property can not be found particularly because models $\mathcal{T}_{t}$ have an ultra-large diameter $D_{t}$ (i.e., $D_{t}\sim|\mathcal{T}_{t}|^{\ln3/\ln5}$) and a greater average shortest path length $\langle\mathcal{W}_{t}\rangle$ (namely, $\langle\mathcal{W}_{t}\rangle\sim|\mathcal{T}_{t}|^{\ln3/\ln5}$) where $|\mathcal{T}_{t}|$ represents vertex number. Next, we determine Pearson correlation coefficient and verify that networks $\mathcal{T}_{t}$ display disassortative mixing structure. In addition, we study random walks on tree networks $\mathcal{T}_{t}$ and derive exact solution to mean hitting time $\langle\mathcal{H}_{t}\rangle$. The results suggest that the analytic formula for quantity $\langle\mathcal{H}_{t}\rangle$ as a function of vertex number $|\mathcal{T}_{t}|$ shows a power-law form, i.e., $\langle\mathcal{H}_{t}\rangle\sim|\mathcal{T}_{t}|^{1+\ln3/\ln5}$. Accordingly, we execute extensive experimental simulations, and demonstrate that empirical analysis is in strong agreement with theoretical results. Lastly, we provide a guide to extend the proposed iterative manner in order to generate more general scale-free tree networks with large diameter. \\


\end{quote}

\vskip 1cm

\section{Introduction}

Nowadays, complex network, as a powerful tool, has proven useful in a great number of applications in different fields ranging from statistic physics, applied mathematic, computer science to chemistry, even to social science \cite{Barabasi-2016}-\cite{Newman-2010}. In particular, the past two decades has seen an upsurge in the area of complex networks. In 1998, Watts and Strogatz established the well-known WS-model for purpose of mimicking small-world property popularly observed on various complex systems \cite{Watts-1998}. Generally speaking, a networked model is considered small-world if it has both smaller diameter and higher clustering coefficient. Almost at the same time, Barab\'{a}si and Albert in 1999 published the famous BA-model in order to capture scale-free feature prevalently encountered in a wide range of complex networks \cite{Barabasi-1999}. A network is called scale-free model when it obeys power-law degree distribution. The core ingredients in BA-model are \emph{growth} and \emph{preferential attachment mechanism}. Indeed, since then, a lot of networked models with small-world and scale-free characteristics \cite{Catanzaro-2005}-\cite{Telesford-2011}, both real-world and synthetic, have been proposed based mainly on thought reflected by two classic models mentioned above. Among which,  Cohen \emph{et al} \cite{Cohen-2003} proved that scale-free networked models are small and, sometimes, ultra-small with respect to relationship between diameter and vertex number. Afterwards, however, some theoretical scale-free models are proved to possess a larger diameter \cite{Zhang-2007}-\cite{Wang-2020-chaos}. As will be shown shortly, our models in this work turn out to follow power-law degree distribution and, at the same time, also have a larger diameter. To put this another way, our scale-free models are not "small" in terms of demonstration in \cite{Cohen-2003}.

As known, there are various kinds of networked models proposed in the last years, such as Apollonian networks \cite{Andrade-2005}, Pseudofractal scale-free web \cite{Dorogovtsev-2002}, Hierarchical networks \cite{Ravasz-2003} and Flower graphs \cite{Diggans-2020}. Tree network, as the simplest connected networked model, has been widely studied in various disciplines \cite{Athreya-2017}-\cite{Ma-2020-T}. For instance, the scale-free BA-tree is immediately obtained by setting the number of edges $m$, which originates from each newly added vertex in BA-model, to $1$ \cite{Barabasi-1999}. Clearly, the scale-free BA-tree is a stochastic model. Similarly, Souza \emph{et al} \cite{Souza-2007} also proposed a family of random recursive trees and found that greedy choice broadens the degree distribution according to the power of choice. In \cite{Gabel-2014}, Gabel \emph{et al} discussed growth tree networks and studied the associated structural parameters. In addition, some trees with deterministic structure have also been studied in detail, such as, geometric fractal growth tree model \cite{Jung-2002}, deterministic uniform recursive tree \cite{Lu-2010}, growing treelike networks \cite{Zhang-2010}, Fibonacci treelike models \cite{Ma-2020-1}, Vicsek fractal \cite{Vicsek-1983}-\cite{Ma-epl-2021} and T-graph \cite{Redner-2001}-\cite{Ma-2020-T}. In this paper, we will introduce a class of new tree networks having interesting features.

Often, in order to better understand the underlying structures on tree networks, some simple yet intriguing dynamics including random walks are considered \cite{Dyer-2020}-\cite{Mokhtar-2013}. In \cite{Vicsek-1983}, Vicsek studied random walks on Vicsek fractal and derived mean hitting time for random walks. Analogously, Zhang \emph{et al} \cite{Zhang-2010} discussed mean first-passage time for random walks on growing treelike networks using spectral technique. More generally, it is challenging to determine exact solution to some structural parameters including mean hitting time on a growth tree network with thousands of vertices. On the other hand, it is of significant interest to derive closed-form solutions to topological parameters of such type. Roughly speaking, there are two main reasons: The first one is that doing so enables one to in depth unveil some interesting structural characteristics behind tree networks under consideration. For example, one can determine whether a given tree network is assortative or not by calculating Pearson correlation coefficient. The other is that it is helpful to develop some effective methods suitable for addressing some specific dynamical problems on networks. Indeed, our discussions in the following several sections will validate the viewpoints above. Here, for instance, when we derive a precise solution to mean hitting time for random walks on tree networks proposed, we do not make use of the commonly used spectral technique but develop a series of simple combinatorial methods. The latter seems more suitable for our tree networks as will be shown later. And, this technique is different from those strategies adopted in most prior work.

The main contributions in this work are as follows:

(1) We introduce an algorithm called Algorithm-I, based on which a family of growth tree networks $\mathcal{T}_{t}$ are iteratively built. By determining degree distribution, we show that networks $\mathcal{T}_{t}$ have scale-free feature. Also, networks $\mathcal{T}_{t}$ turn out to have a great diameter. In addition, an analogy, scale-free tree networks $\mathcal{T}^{\star}_{t}$ that have the vertex number with networks $\mathcal{T}_{t}$ in common, is also constructed and discussed. These two networks are proved to display disassortative mixing structure. 

(2) We study random walks on scale-free networks $\mathcal{T}_{t}$, and derive the exact solution of mean hitting time via proposing a series of combinatorial manners. Also, the solution of mean hitting time on networks $\mathcal{T}^{\star}_{t}$ is obtained analytically. The results show that diameter plays a key role in determination of mean hitting time of networks $\mathcal{T}_{t}$ and $\mathcal{T}^{\star}_{t}$.

(3) We conduct extensive computer simulations, which confirms that empirical analysis is perfectly consistent with the theoretical results. At the same time, we also provide some guidelines for generating tree networks that can exhibit other structural properties when considering random walks, including scale-free tree network that nearly achieves theoretical upper bound for mean hitting time.

The rest of this paper is organized as follows. First, we introduce some terminologies in Section 2. Section 3 aims at presenting the concrete procedure for generating our tree network $\mathcal{T}_{t}$. At the same time, the number of vertices of network $\mathcal{T}_{t}$ is derived in an iterative fashion. Next, we estimate several significant structural parameters, i.e., degree distribution, diameter and Pearson correlation coefficient, in Section 4. The results suggest that (1) tree network $\mathcal{T}_{t}$ follows power-law degree with exponent larger than $3$, (2) "small-world" property is not observed on network $\mathcal{T}_{t}$, and (3) network $\mathcal{T}_{t}$ is disassortative mixing. In Section 5, we study random walks on tree network $\mathcal{T}_{t}$, and derive the closed-form solution of mean hitting time. The corresponding analytic formula for mean hitting time as a function of vertex number $|\mathcal{T}_{t}|$ exhibits a power-law form. Note that throughout the preceding two sections, we also in depth discuss the adjoint tree model of tree network $\mathcal{T}_{t}$ which is viewed as $\mathcal{T}^{\star}_{t}$. Related work is reviewed briefly in Section 6. Finally, we close this work by drawing the conclusion in Section 7.

\section{Terminologies}

This section aims at introducing some fundamental terminologies used in this paper, including degree distribution, fractal and fractal dimension. Accordingly, some notations are collected in table 1.

\subsection{Degree distribution }

Let $\mathcal{G}(\mathcal{V},\mathcal{E})$ denote a graph that is an ordered pair ($\mathcal{V},\mathcal{E}$) consisting of a set $\mathcal{V}$ of vertices and a set $\mathcal{E}$ of edges running between vertices \cite{Bondy-2008}. Accordingly, we make use of symbols $|\mathcal{V}|$ and $|\mathcal{E}|$ to represent vertex number and edge number, respectively. If $\mathcal{N}_{k}$ is referred to as the total number of vertices with degree $k$, then the probability of selecting a degree $k$ vertex at random from graph $\mathcal{G}(\mathcal{V},\mathcal{E})$, denoted by $P(k)$, is given by

\begin{equation}\label{eqa:MF-TKDE-2-2-1}
P(k)=\frac{\mathcal{N}_{k}}{|\mathcal{V}|}.
\end{equation}
More generally, due to a fact that degree values in network are discrete in form, one is concentrated with cumulative degree distribution $P_{cum}(k)$ of network in question which is as follows
\begin{equation}\label{eqa:MF-TKDE-2-2-2}
P_{cum}(k)=\frac{\sum_{k_{i}\geq k}\mathcal{N}_{k_{i}}}{|\mathcal{V}|}.
\end{equation}
Clearly, with simple arithmetics, we can derive a relation between quantities $P(k)$ and $P_{cum}(k)$ as below

\begin{equation}\label{eqa:MF-TKDE-2-2-3}
P(k)=\dfrac{d}{dk}P_{cum}(k).
\end{equation}

\subsection{Fractal and fractal dimension }

In fact, there is a long history in the study of fractal. For example, a fractal has been defined in Mandelbrot (1982) as a rough or fragmented geometric shape that can be split into parts, each of which is (at least approximately) a reducedsize copy of the whole \cite{Oliver-2013}. In the context of fractal, one fundamental yet important parameter for characterizing a fractal is the fractal dimension, denoted by $d_{f}$. It is a measure of the extent to which a structure exceeds its base dimension to fill the next dimension. In general, there are many distinct manners used to evaluate fractal dimension $d_{f}$ for a given fractal. In this paper, we make use of the following definition: Consider a perfectly self-similar fractal with $r$ similar pieces, each of which is scaled down by a factor $s$. And then, the fractal dimension for the fractal under consideration is viewed as

\begin{equation}\label{eqa:MF-TKDE-2-1-1}
d_{f}=\frac{\ln r}{\ln s}.
\end{equation}

To sum up, before beginning with discussions, we list out some notations used later for reader to easy understand theoretical analysis in the rest of this paper. These notations are shown in table 1.

\begin{table}                                                                                                                                                                                                                                                                                          
\centering                                                                                                                                                                                                                                                                                             
\caption{ Summary of notations.}                                                                                                                                                                                                         
\label{table}                                                                                                                                                                                                                                                                                          
\begin{tabular}{|c|c|}
\hline                                                                                                                                                                                                                                                                                                 
  Notation&  Meaning  \\                                                                                                                                                                                                                                                                                
\hline                                                                                                                                                                                                                                                                                                 
  $\mathcal{G}(\mathcal{V},\mathcal{E})$&  Graph (or Network), vertex set $\mathcal{V}$, edge set $\mathcal{E}$  \\                                                                                                                                                                                                                                                           
\hline                                                                                                                                                                                                                                                                                                 
$\mathcal{T}$&   Tree \\                                                                                                                                                                                                                                                                               
\hline                                                                                                                                                                                                                                                                                                 
  $\mathcal{P}_{uv}$&   A path connecting vertex $u$ to vertex $v$\\                                                                                                                                                                                                                                                                        
\hline                                                                                                                                                                                                                                                                                                 
$e_{uv}$(or $uv$)&  An edge connecting vertices $u$ and $v$  \\                                                                                                                                                                                                                                                            
\hline                                                                                                                                                                                                                                                                                                 
  $d_{uv}$&   Distance between vertices $u$ and $v$\\ 
\hline                                                                                                                                                                                                                                                                                                 
$D$&  Diameter of network \\                                                                                                                                                                                                                                                            
\hline                                                                                                                                                                                                                                                                                                 
  $d_{f}$&   Fractal dimension of fractal network\\                                                                                                                                                                                                                                                                                                                                                                                                                                                                                                                           
\hline                                                                                                                                                                                                                                                                                                 
$r$&  Pearson correlation coefficient\\                                                                                                                                                                                                                                                                                                                                                                                                                                                                                                                             
\hline  
$\mathcal{N}_{k}$&  The number of vertices with degree $k$\\                                                                                                                                                                                                                                                                                                                                                                                                                                                                                                                             
\hline 
$\Omega_{u}$&  The neighboring set $\Omega_{u}$ of vertex $u$ in graph\\                                                                                                                                                                                                                                                                                                                                                                                                                                                                                                                             
\hline   
$A_{v}(t)$&  Set of new vertices inserted on all the existing edges incident with vertex $v$ in $\mathcal{T}_{t}$ \\                                                                                                                                                                                                                                                                                                                                                                                                                                                                                                                             
\hline  
$B_{v}(t)$&  Set of new vertices as leaves directly connected to vertex $v$ in tree $\mathcal{T}_{t}$\\                                                                                                                                                                                                                                                                                                                                                                                                                                                                                                                             
\hline 
$k_{v}(t)$&  Degree of vertex $v$ in tree $\mathcal{T}_{t}$\\                                                                                                                                                                                                                                                                                                                                                                                                                                                                                                                             
\hline                                                                                                                                                                                                                                                                                              
$P(k)$&  Degree distribution of network\\                                                                                                                                                                                                                                                     
\hline 
$P_{cum}(k)$&  Cumulative degree distribution of network \\                                                                                                                                                                                                                                                            
\hline                                                                                                                                                                                                                                                                                                 
$\mathcal{H}_{u\rightarrow v}$&  Hitting time from vertex $u$ to $v$ for random walks on network\\                                                                                                                                                                                                                                                                                                                                                                                                                                                                                                                           
\hline                                                                                                                                                                                                                                                                                                 
$\langle\mathcal{H}_{\mathcal{G}}\rangle$&  Mean hitting time for random walks on network $\mathcal{G}(\mathcal{V},\mathcal{E})$\\                                                                                                                                                                                                                                                                                                                                                                                                                                                                                                                             
\hline                                                                                                                                                                                                                                                                                                 
$\langle\mathcal{W}_{\mathcal{G}}\rangle$&  Average shortest path length of network $\mathcal{G}(\mathcal{V},\mathcal{E})$\\                                                                                                                                                                                                                                                     
\hline                                                                                                                                                                                                                                                                                                  
\end{tabular}                                                                                                                                                                                                                                                                                          
\end{table}

\section{Construction}

In this section, we propose an algorithm (called Algorithm-I for convenience) for creating tree network $\mathcal{T}_{t}$, and then study some fundamental parameters. Part of them are helpful to analyze other structural quantities as will be shown shortly. At the meantime, the corresponding adjoint tree network $\mathcal{T}^{\star}_{t}$ is also established. Note that the exact solutions of some structural parameters of tree $\mathcal{T}^{\star}_{t}$ are given but the associated detailed derivations are omitted. This is because they can be conveniently obtained in a similar manner to that for tree network $\mathcal{T}_{t}$.

\subsection{Algorithm-I}

First of all, it should be mentioned that an arbitrary tree $\mathcal{T}$ can be selected as input to the algorithm proposed below. Nonetheless, as tried in the rich literature \cite{Zhang-2010,Vicsek-1983,Agliari-2008}, we only consider a specific case, i.e., a single edge serving as the seed of Algorithm-I, when studying degree distribution. In other cases, we still consider general setting, namely, an arbitrary tree as seed.

At $t=0$, the seminal model is an arbitrary tree $\mathcal{T}$, referred to as $\mathcal{T}_{0}$. For our purpose, we denote by $|\mathcal{T}_{0}|$ vertex number.

At $t=1$, the second tree $\mathcal{T}_{1}$ is obtained from preceding tree $\mathcal{T}_{0}$ by performing the following two operations:

\begin{itemize}
\item We first insert one new vertex into each edge incident with each vertex $v$ belonging to tree $\mathcal{T}_{0}$.

\item And then, we connect $k_{v}(0)$ new vertices as leaves to each vertex $v$ of tree $\mathcal{T}_{0}$. Here, symbol $k_{v}(0)$ represents the degree of vertex $v$ in tree $\mathcal{T}_{0}$.

\emph{Remark}---Such a manipulation is easy to extend to many more complicated operations. Yet, our goal is to introduce this idea and thus some more intriguing extensions are left as investigations in future research.

\end{itemize}

At $t\geq2$, tree $\mathcal{T}_{t}$ is generated based on the previous tree $\mathcal{T}_{t-1}$ using the same operation as mentioned above.

To make our description more concrete, an illustrative example is shown in panel (a) of Fig.1 where a single edge is chose as a seed of Algorithm-I.

\begin{figure}
\centering
  \includegraphics[height=6cm]{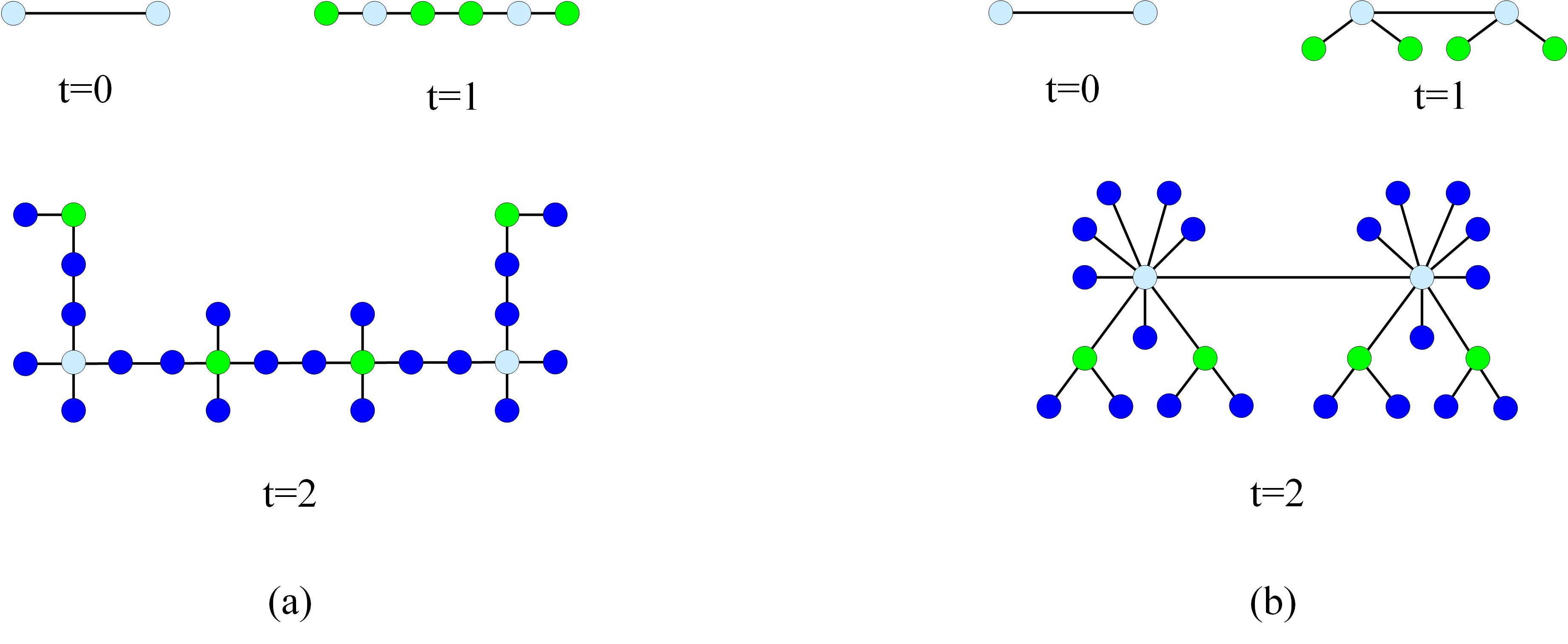}\\
\caption{(Color online) The diagram of tree networks $\mathcal{T}_{t}$ and $\mathcal{T}^{\star}_{t}$. The first three generations of tree $\mathcal{T}_{t}$ are plotted in panel (a). Among of them, we make use of distinct colors to highlight vertices added at different time steps. Similarly, we also show the first three generations of tree $\mathcal{T}^{\star}_{t}$ in panel (b). }
\end{figure}

\subsection{Vertex number}

Clearly, our tree network $\mathcal{T}_{t}$ is also built upon an iterative method, a generative approach that has been frequently used to create various networked models in the past years \cite{Ma-2018,Zhang-2007,Dorogovtsev-2002,Ravasz-2003}. Accordingly, this kind of constructive way enables calculations of some fundamental topological parameters on network created, such as, vertex number. Below, we provide a simple derivation of vertex number $|\mathcal{T}_{t}|$ by virtue of the concrete growth way of tree network $\mathcal{T}_{t}$.

\textbf{Theorem 3.1} The exact solution of vertex number $|\mathcal{T}_{t}|$ is given by

\begin{equation}\label{eqa:MF-1-1}
|\mathcal{T}_{t}|=5^{t}(|\mathcal{T}_{0}|-1)+1.
\end{equation}

\emph{Proof} From Algorithm-I, we can see that in fact, there are $2k_{v}(t-1)$ new vertices created for each vertex $v$ of tree $\mathcal{T}_{t-1}$ at time step $t$($>0$). This leads to the following relation

$$|\mathcal{T}_{t}|=\sum_{v\in\mathcal{T}_{t-1}}(1+2k_{v}(t-1)).$$

With a trivial fact that the summation of degrees of all vertices in a tree $\mathcal{T}$ is equal to two times (vertex number $|\mathcal{T}|$ minus one), the exact solution to vertex number $|\mathcal{T}_{t}|$ is obtained as shown in Eq.(\ref{eqa:MF-1-1}). \qed

Besides that, we want to make more detailed demonstration about Algorithm-I in order to enable the next calculations. Specifically, there are two classes of new vertices added at each time step. As a case study, it is easy to see that $k_{v}(t-1)$ new vertices are inserted on all the existing edges incident with vertex $v$ of $\mathcal{T}_{t-1}$ with each being in each edge. Thus, each of vertices of this type has initial degree $2$. On the contrary, the other $k_{v}(t-1)$ new vertices are directly connected to vertex $v$ itself, which means that each of them has initial degree $1$. For brevity, the former vertices associated with vertex $v$ are grouped into set $A_{v}(t)$, and, similarly, notation $B_{v}(t)$ is used to indicate a set of the latter vertices. To put this another way, at time step $t$, tree network $\mathcal{T}_{t}$ has three types of vertex sets, i.e., $\mathcal{T}_{t-1}$, $\bigcup_{v\in\mathcal{T}_{t-1}}A_{v}(t)$ and $\bigcup_{v\in\mathcal{T}_{t-1}}B_{v}(t)$. For brevity, we define $A(t):=\bigcup_{v\in\mathcal{T}_{t-1}}A_{v}(t)$ and $B(t)=\bigcup_{v\in\mathcal{T}_{t-1}}B_{v}(t)$.

So far, we finish the construction of tree network $\mathcal{T}_{t}$. To make further progress, if we adapt Algorithm-I to let $k_{v}(t-1)$ new vertices in set $A_{v}(t)$ not be placed on edges incident with vertex $v$ but also connect to vertex $v$ itself directly at each time step $t$, then we obtain another tree network, called $\mathcal{T}^{\star}_{t}$ for our purpose. Panel (b) of Fig.1 shows the first three generations of tree network $\mathcal{T}^{\star}_{t}$. It is obvious to understand that tree network $\mathcal{T}^{\star}_{t}$ shares an identical vertex number with tree network $\mathcal{T}_{t}$ when a common tree is selected as their seeds, which is shown in the following theorem.

\textbf{Theorem 3.2} The exact solution of vertex number $|\mathcal{T}^{\star}_{t}|$ is given by

\begin{equation}\label{eqa:MF-1-2}
|\mathcal{T}^{\star}_{t}|=5^{t}(|\mathcal{T}_{0}|-1)+1.
\end{equation}

\emph{Proof} This is verified in the same method as that used in the proof of theorem 3.1 and we thus omit it here. \qed

Nonetheless, both tree networks $\mathcal{T}_{t}$ and $\mathcal{T}^{\star}_{t}$ still exhibit significant difference on underlying structure when we discuss some other structural parameters as will be shown shortly. It is worth noting that a more general version regarding to tree network $\mathcal{T}^{\star}_{t}$ has been studied in detail in \cite{Jung-2002}.

\section{Several significant topological parameters}

From now on, we will focus on three significant structural parameters, degree distribution, diameter and Pearson correlation coefficient, on tree network $\mathcal{T}_{t}$. One of important reasons for this is that these parameters have been widely used to study the topological structure of various networked models in the literature \cite{Barabasi-2016}-\cite{Barabasi-1999}. For example, it is using degree distribution that Barab\'{a}si \emph{et al} revealed the scale-free feature of a wide range of networks \cite{Barabasi-1999}. With diameter, a great number of networks turn out to be small-world \cite{Watts-1998}. Many biological and technological networks have been proven to belong to the disassortative mixing family by determining Pearson correlation coefficient \cite{Newman-2002}. More generally, determining structural parameters on networked models is useful to understand the underlying structure on networks. As will see in the following, our tree network $\mathcal{T}_{t}$ indeed shows some characteristics associated with both degree distribution and diameter that can not be found in some previous tree models \cite{Lu-2010}-\cite{Ma-2020-T}.

\subsection{Degree distribution}

Degree distribution, as an important index, has been popularly adopted to determine whether a given network is scale-free or not in the study of complex networks \cite{Barabasi-2016}. Informally, one always uses $P(k)$ to define the probability of choosing at random uniformly a degree $k$ vertex in a given network $\mathcal{G}(\mathcal{V},\mathcal{E})$. In the thermodynamic limit, probability $P(k)$ is referred to as the ratio of the total number of degree $k$ vertices $\mathcal{N}_{k}$ and vertex number $|\mathcal{V}|$. According to the fact that vertex degrees of network $\mathcal{G}(\mathcal{V},\mathcal{E})$ are discrete values in form, the cumulative degree distribution $P_{cum}(k)$ is often exploited, which is by definition written as

$$P_{cum}(k)=\frac{\sum_{k'\geq k}\mathcal{N}_{k'}}{|\mathcal{V}|}.$$

\textbf{Theorem 4.1} Tree network $\mathcal{T}_{t}$ follows power-law degree distribution as follows

\begin{equation}\label{eqa:MF-2-1-0}
P(k)\sim k^{-(1+\gamma)}, \quad \gamma=\frac{\ln5}{\ln2}.
\end{equation}

\emph{Proof} Clearly, it suffices to rank all the vertices of tree network $\mathcal{T}_{t}$ with respect to vertex degree in order to derive cumulative degree distribution $P_{cum}(k)$. Now, let us keep this in mind and then recall the concrete procedure in Algorithm-I. For one vertex $v$ added into tree $\mathcal{T}_{t}$ at time step $t_{i}$, we can see that its degree $k_{v}(t)$ might take two distinct values. If vertex $v$ is in set $A(t_{i})$, then the corresponding degree $k_{v}(t)$ is equal to $2^{t-t_{i}+1}$. On the other hand, degree $k_{v}(t)$ is equivalent to $2^{t-t_{i}}$ when vertex $v$ belongs to set $B(t_{i})$. Based on this, we can obtain the following expression

\begin{equation}\label{eqa:MF-2-1-1}
P_{cum}(k)=\frac{2+\sum_{j=1}^{t_{i}}|A(j)|+\sum_{j=1}^{t_{i}-1}|B(j)|}{5^{t}+1}
\end{equation}
where $|A|$ indicates the cardinality of a set $A$ and we have made use of a single edge as seed as stated above. Substituting $k=2^{t-t_{i}+1}$ into Eq.(\ref{eqa:MF-2-1-1}) yields

$$P_{cum}(k)\sim k^{-\gamma},\quad \gamma=\frac{\ln5}{\ln2}.$$

Performing derivative on both hand sides in the above expression produces

\begin{equation}\label{eqa:MF-2-1-2}
P(k)\sim k^{-(1+\gamma)},
\end{equation}
which suggests that our tree $\mathcal{T}_{t}$ follows power-law degree distribution. In other words, scale-free feature is observed on tree $\mathcal{T}_{t}$ in the limit of large graph size. This completes the proof of theorem 4.1. \qed

It is worth noting that power-law exponent $1+\gamma$ is certainly larger than $3$. This is a phenomenon that is seldom reported in published tree networks with scale-free feature. In essence, it is also confirmed in terms of "preferential attachment" in the following observation. In addition, we prove by definition that tree network $\mathcal{T}_{t}$ is in fact fractal model, and, accordingly, obtain that its fractal dimension $d_{f}$ is equal to $\ln5/\ln3$.   

Analogously, we can without difficulty derive degree distribution of tree $\mathcal{T}^{\star}_{t}$, which is shown in the coming theorem.

\textbf{Theorem 4.2} Tree network $\mathcal{T}^{\star}_{t}$ follows power-law degree distribution as follows

\begin{equation}\label{eqa:MF-2-1-3}
P^{\star}(k)\sim k^{-(1+\gamma^{\star})}, \quad \gamma^{\star}=\frac{\ln5}{\ln3}.
\end{equation}

\emph{Proof} Due to a similar analysis as above, we omit detailed proof. \qed

From which we can see that tree $\mathcal{T}^{\star}_{t}$ obeys power-law degree distribution as well. However, its power-law exponent is smaller than $3$. A similar result has been reported in \cite{Jung-2002}.

Apparently, the above analysis of degree distribution implies that tree $\mathcal{T}^{\star}_{t}$ has a different underlying structure from tree $\mathcal{T}_{t}$. Yet, it should be mentioned that in the large graph size limit, they have an identical average degree regardless of power-law exponent. In addition, the both types of generative ways to create scale-free tree networks $\mathcal{T}_{t}$ and $\mathcal{T}^{\star}_{t}$ are within the category of "preferential attachment mechanism" introduced by Barab\'{a}si \emph{et al} \cite{Barabasi-1999}. In the past two decades, most papers focus on the generative approach to generating tree $\mathcal{T}^{\star}_{t}$ \cite{Ma-2020-1,Kim-2004}. As a result, almost all published tree networks with scale-free feature are always proved to have an exponent no more than $3$. Therefore, this work provides an effective way to build up scale-free tree networks with intriguing features, such as power-law exponent larger than $3$. Additionally, we provide some deep discussions about the proposed generative way in the end of Section 5.

\subsection{Diameter}

In the jargon of graph theory, diameter of network $\mathcal{G}(\mathcal{V},\mathcal{E})$, denoted by $D$, is the distance of longest shortest path between all possible pairs of vertices. As one of most fundamental topological parameters, diameter has been proved useful in a variety of applications \cite{Dai-2019,Bianchin-2015}. For example, it is in general utilized to measure information delay in communication network. As known, small-world phenomena on complex networks are certified by means of diameter \cite{Watts-1998}. Next, let us pay attention to discussion about diameter of tree network $\mathcal{T}_{t}$.

\textbf{Theorem 4.3} The closed-form solution of diameter of tree network $\mathcal{T}_{t}$ is given by
\begin{equation}\label{eqa:MF-2-2-0}
D_{t}=3^{t}(D_{0}+1)-1.
\end{equation}

\emph{Proof} As shown in Algorithm-I, two new vertices are inserted into each pre-existing edge at each time step. This means that diameter $D_{t}$ of tree $\mathcal{T}_{t}$ satisfies the following relation, namely,

$$D_{t}=2+3D_{t-1}.$$

With initial condition $D_{0}$ that is diameter of the seed $\mathcal{T}_{0}$, it is straightforward to write

\begin{equation}\label{eqa:MF-2-2-1}
D_{t}=3^{t}(D_{0}+1)-1,
\end{equation}
which is completely the same as Eq.(\ref{eqa:MF-2-2-0}). We complete the proof of theorem 4.3. \qed

Obviously, unlike statement in \cite{Cohen-2003}, diameter $D_{t}$ of tree network $\mathcal{T}_{t}$ is much more than $\ln|\mathcal{T}_{t}|$ and follows $D_{t}\sim|\mathcal{T}_{t}|^{\ln3/\ln5}$. This means that our tree network $\mathcal{T}_{t}$ has no small-world property. As is known to us all, diameter is a coarser index for deciding small-world property. A finer one is average shortest path length. One thing to note is that in the end of Section 5, we derive the closed-form solution of average shortest path length on tree network $\mathcal{T}_{t}$ (which is shown in Eq.(\ref{eqa:MF-A-1})), which certainly suggests that our network is indeed not small-world.

As above, we also obtain exact solution of diameter $D^{\star}_{t}$ of tree network $\mathcal{T}^{\star}_{t}$ in an iterative fashion and write it in the following form

\begin{equation}\label{eqa:MF-2-2-2}
D^{\star}_{t}=D_{0}+2t,
\end{equation}
in which $D_{0}$ has the same meaning as mentioned above. It is clear to see that tree $\mathcal{T}^{\star}_{t}$ meets the criteria in \cite{Cohen-2003}, i.e., $D^{\star}_{t}\sim\ln|\mathcal{T}^{\star}_{t}|$. A similar result is shown in \cite{Jung-2002}. In some published papers \cite{Zhang-2010}, trees of such type are considered as models with small-world property. Strictly speaking, an arbitrary tree network can not be small-world mainly due to the corresponding concept, i.e., small diameter and higher clustering coefficient. Tree network always has zero clustering coefficient.

From Eqs.(\ref{eqa:MF-2-2-1}) and (\ref{eqa:MF-2-2-2}), we can also state that both types of generative ways introduced in Section 2 have significantly different impacts on structure of tree networks. Intuitively, the way to create tree $\mathcal{T}_{t}$ enlarges the underlying structure of network more remarkably than that for creating tree $\mathcal{T}^{\star}_{t}$. Therefore, it is a useful manner in which one may build up some networked models having a quite large diameter. Note also that a similar operation, which makes the growth of diameter of network under consideration more rapid, has been used to construct many other networked models of great interest \cite{Zhang-2007}-\cite{Wang-2020-chaos}, including scale-free networks with both density feature and large diameter \cite{Ma-2020-PRE} and fractal and transfractal recursive scale-free nets \cite{Rozenfeld-2007}. While the proposed models in \cite{Ma-2020-PRE,Rozenfeld-2007} display some interesting properties, such as, scale-free feature, they are not treelike. In this paper, however, our focus is on tree networks. Generally speaking, this work can expend our understanding about deterministic networked models, especially, treelike networks.

\subsection{Pearson correlation coefficient}

As the final structural parameter discussed in this section, Pearson correlation coefficient $r$ of network $\mathcal{G}=(\mathcal{V},\mathcal{E})$, firstly studied by Newman in \cite{Newman-2002}, is defined as
\begin{equation}\label{eqa:MF-2-3-0}
r=\frac{|\mathcal{E}|^{-1}\sum\limits_{e_{ij}\in \mathcal{E}} k_{i}k_{j}-\left[|\mathcal{E}|^{-1}\sum\limits_{e_{ij}\in \mathcal{E}} \frac{1}{2}(k_{i}+k_{j})\right]^{2}}{|\mathcal{E}|^{-1}\sum\limits_{e_{ij}\in \mathcal{E}} \frac{1}{2}(k^{2}_{i}+k^{2}_{j})-\left[|\mathcal{E}|^{-1}\sum\limits_{e_{ij}\in \mathcal{E}} \frac{1}{2}(k_{i}+k_{j})\right]^{2}},
\end{equation}
in which $k_{i}$ is degree of vertex $i$ and $e_{ij}$ denotes an edge connecting vertex $i$ to $j$. 

By definition, one can clearly understand that parameter $r$ in essence measures tendency of connections taking place between vertices of network $\mathcal{G}=(\mathcal{V},\mathcal{E})$. Case of $r>0$ means that there is a clear preference for vertex to link to other vertex like itself with regard of vertex degree. And, case of $r<0$ indicates the opposite consequence. $r=0$ is the critical value. Now, let us calculate this parameter on tree network $\mathcal{T}_{t}$ and decide what kind of case tree network $\mathcal{T}_{t}$ belongs to.

\textbf{Theorem 4.4} The exact solution of Pearson correlation coefficient $r_{t}$ of tree network $\mathcal{T}_{t}$ is given by

\begin{equation}\label{eqa:MF-2-3-1}
r_{t}=\frac{\frac{\Gamma_{t}(1)}{5^{t}}-\left[\frac{\Gamma_{t}(2)}{2\times 5^{t}}\right]^{2}}
{\frac{\Gamma_{t}(1)}{2\times 5^{t}}-\left[\frac{\Gamma_{t}(2)}{2\times 5^{t}}\right]^{2}},
\end{equation}
here

\begin{subequations}
\label{eq:whole}
\begin{eqnarray}
\begin{aligned}\Gamma_{t}(1)=6\times\left[4\times2^{t-1}+\sum_{i=1}^{t-2}2^{i+1}(5^{t-2-i}+5^{t-1-i})+2\times 5^{t-1}\right]+4\times5^{t-1},
\end{aligned}\label{subeq:MF-1-1}
\end{eqnarray}
\begin{equation}
\begin{aligned}\Gamma_{t}(2)=4\times\left[4\times2^{t-1}+\sum_{i=1}^{t-2}2^{i+1}(5^{t-2-i}+5^{t-1-i})+2\times 5^{t-1}\right]+2\times5^{t},
\end{aligned}\label{subeq:MF-1-2}
\end{equation}
\begin{equation}
\begin{aligned}\Gamma_{t}(3)=8\times\left[4\times2^{t-1}+\sum_{i=1}^{t-2}2^{i+1}(5^{t-2-i}+5^{t-1-i})+2\times 5^{t-1}\right]+18\times5^{t-1}.
\end{aligned}\label{subeq:MF-1-3}
\end{equation}
\end{subequations}
Note that we have employed the assumption that a single edge is selected as seed for creating tree network $\mathcal{T}_{t}$. 

\emph{Proof} First of all, we introduce three notations as follows 

$$\Gamma_{t}(1):=\sum_{e_{uv}\in \mathcal{T}_{t}}k_{u}(t)\times k_{v}(t),\qquad \Gamma_{t}(2):=\sum_{e_{uv}\in \mathcal{T}_{t}}\left(k_{u}(t)+k_{v}(t)\right), \qquad \Gamma_{t}(3):=\sum_{e_{uv}\in \mathcal{T}_{t}}\left([k_{u}(t)]^{2}+[k_{v}(t)]^{2}\right).$$

Based on the concrete construction of tree network $\mathcal{T}_{t}$, we obtain 

\begin{equation}\label{eqa:MF-2-3-2}
\begin{aligned}\Gamma_{1}(t)&=\sum_{u\in \mathcal{T}_{t-1}}2k_{u}(t-1)\times2k_{u}(t-1)+\sum_{u\in \mathcal{T}_{t-1}}2k_{u}(t-1)\times k_{u}(t-1)+\sum_{e_{uv}\in \mathcal{T}_{t-1}}4\\
&=6\sum_{u\in \mathcal{T}_{t-1}}[k_{u}(t-1)]^{2}+4\left(|\mathcal{T}_{t-1}|-1\right)\\
&=6\times\left[4\times2^{t-1}+\sum_{i=1}^{t-2}2^{i+1}(5^{t-2-i}+5^{t-1-i})+2\times 5^{t-1}\right]+4\times5^{t-1},
\end{aligned}
\end{equation}
in which we have used Eq.(\ref{eqa:MF-1-1}). 

Analogously, the following two expressions are easy to check.

\begin{equation}\label{eqa:MF-2-3-3}
\begin{aligned}\Gamma_{2}(t)&=\sum_{u\in \mathcal{T}_{t-1}}[2k_{u}(t-1)+2]\times k_{u}(t-1)+\sum_{u\in \mathcal{T}_{t-1}}[2k_{u}(t-1)+1]\times k_{u}(t-1)+\sum_{e_{uv}\in \mathcal{T}_{t-1}}4\\
&=4\sum_{u\in \mathcal{T}_{t-1}}[k_{u}(t-1)]^{2}+10\left(|\mathcal{T}_{t-1}|-1\right)\\
&=4\times\left[4\times2^{t-1}+\sum_{i=1}^{t-2}2^{i+1}(5^{t-2-i}+5^{t-1-i})+2\times 5^{t-1}\right]+2\times5^{t},
\end{aligned}
\end{equation}

\begin{equation}\label{eqa:MF-2-3-4}
\begin{aligned}\Gamma_{3}(t)&=\sum_{u\in \mathcal{T}_{t-1}}\left(4[k_{u}(t-1)]^{2}+4\right)\times k_{u}(t-1)+\sum_{u\in \mathcal{T}_{t-1}}\left(4[k_{u}(t-1)]^{2}+1\right)\times k_{u}(t-1)+\sum_{e_{uv}\in \mathcal{T}_{t-1}}8\\
&=8\sum_{u\in \mathcal{T}_{t-1}}[k_{u}(t-1)]^{2}+18\left(|\mathcal{T}_{t-1}|-1\right)\\
&=8\times\left[4\times2^{t-1}+\sum_{i=1}^{t-2}2^{i+1}(5^{t-2-i}+5^{t-1-i})+2\times 5^{t-1}\right]+18\times5^{t-1}.
\end{aligned}
\end{equation}

Lastly, armed with Eq.(\ref{eqa:MF-1-1}) and Eqs.(\ref{eqa:MF-2-3-2})-(\ref{eqa:MF-2-3-4}), we obtain the same expression as shown in theorem 4.3. This completes the proof of Eq.(\ref{eqa:MF-2-3-1}). \qed 

Similarly, the closed-form solution of Pearson correlation coefficient $r^{\star}_{t}$ of tree network $\mathcal{T}^{\star}_{t}$ is easily derived and shown in the coming form.

\begin{equation}\label{eqa:MF-2-3-5}
r^{\star}_{t}=\frac{\frac{\Gamma^{\star}_{t}(1)}{5^{t}}-\left[\frac{\Gamma^{\star}_{t}(2)}{2\times 5^{t}}\right]^{2}}
{\frac{\Gamma^{\star}_{t}(1)}{2\times 5^{t}}-\left[\frac{\Gamma^{\star}_{t}(2)}{2\times 5^{t}}\right]^{2}},
\end{equation}
here

\begin{subequations}
\label{eq:whole}
\begin{eqnarray}
\begin{aligned}\Gamma^{\star}_{t}(1)=9^{t}+6\sum_{i=0}^{t-1}9^{i}\left(2\times3^{2(t-1-i)}+4\sum_{j=0}^{t-2-i}5^{j}\times3^{2(t-2-j)}\right),
\end{aligned}\label{subeq:MF-2-1}
\end{eqnarray}
\begin{equation}
\begin{aligned}\Gamma^{\star}_{t}(2)=2\times3^{t}+6\sum_{i=0}^{t-1}3^{i}\left(2\times3^{2(t-1-i)}+4\sum_{j=0}^{t-2-i}5^{j}\times3^{2(t-2-j)}\right)+4\sum_{i=0}^{t-1}3^{i}\times5^{t-1-i},
\end{aligned}\label{subeq:MF-2-2}
\end{equation}
\begin{equation}
\begin{aligned}\Gamma^{\star}_{t}(3)=2\times9^{t}+18\sum_{i=0}^{t-1}9^{i}\left(2\times3^{3(t-1-i)}+4\sum_{j=0}^{t-2-i}5^{j}\times3^{3(t-2-j)}\right)+4\sum_{i=0}^{t-1}9^{i}\times5^{t-1-i}.
\end{aligned}\label{subeq:MF-2-3}
\end{equation}
\end{subequations}

As above, we omit the detailed derivation of parameter $r^{\star}_{t}$. Yet, to show the scaling of Pearson correlation coefficients $r_{t}$ and $r^{\star}_{t}$ in the large graph size limit, we feed networks $\mathcal{T}_{t}$ and $\mathcal{T}^{\star}_{t}$ into computer and then gain an illustrative outline in Fig.2. It is clear to see that parameters $r_{t}$ and $r^{\star}_{t}$ are constantly negative and approach $0$ in the limit of large-$t$, which implies that two tree networks $\mathcal{T}_{t}$ and $\mathcal{T}^{\star}_{t}$ possess disassortative structure. In another word, two types of tree networks are disassortative. It should be mentioned that except for time step $t=3$, $r_{t}$ is always smaller than $r^{\star}_{t}$, which suggests that the generative manner for tree network $\mathcal{T}_{t}$ is beneficial to disassortative structure of growing tree network.

\begin{figure}
\centering
  \includegraphics[height=6cm]{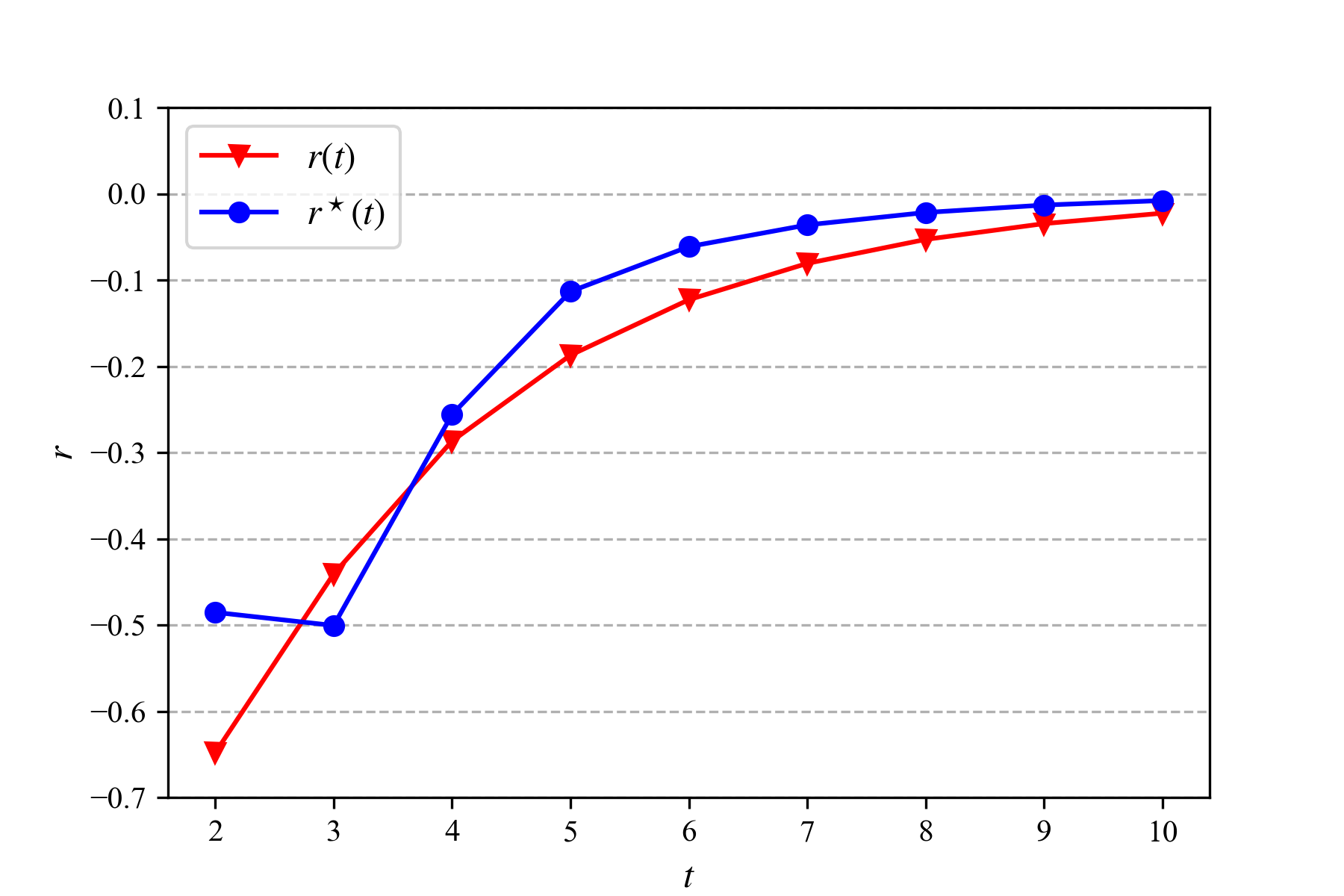}\\
\caption{(Color online) The diagram of Pearson correlation coefficients $r_{t}$ and $r^{\star}_{t}$ where a single edge serves as the seed. The straight line indicates analytical solutions for quantities $r_{t}$ and $r^{\star}_{t}$. Triangle and circle represents, respectively, the numerical value of $r_{t}$ and $r^{\star}_{t}$ produced by computer simulation. Clearly, two parameters are constantly negative and approach $0$ in the limit of large-$t$. }
\end{figure}

Until now, we study some significant structural parameters on tree network $\mathcal{T}_{t}$ in more detail. The results demonstrate that (1) the proposed model $\mathcal{T}_{t}$ follows power-law degree distribution with exponent more than $3$, and (2) model $\mathcal{T}_{t}$ has a larger diameter and thus has no small-world property. At the same time, model $\mathcal{T}_{t}$ possesses disassortative structure due to Pearson correlation coefficient $r_{t}$ smaller than zero. In the next section, we will discuss random walks on tree network $\mathcal{T}_{t}$ in order to further probe the topological structure of network $\mathcal{T}_{t}$.

\section{Random walks on tree network }

Random walks, as the simplest representative of discrete-time unbiased diffusion phenomena, have received increasing attention from various fields ranging from applied mathematics, statistic physics, computer science to bio-chemistry \cite{Redner-2001,Noh-2004,Beveridge-2013,MF-2020,Oliver-2013,Hwang-2012,Perkins-2014,Nikolakopoulos-2020}, and so forth. A walker performing random walks on network $\mathcal{G}(\mathcal{V},\mathcal{E})$ moves with a probability $P_{u\rightarrow v}=1/k_{u}$ from his current position $u$ to each candidate $v$ in the neighboring set $\Omega_{u}$ in a step. Such a dynamic process is in general encoded by the corresponding transition matrix $\mathbf{P}_{\mathcal{G}}=\mathbf{D}_{\mathcal{G}}^{-1}\mathbf{A}_{\mathcal{G}}$ of network $\mathcal{G}(\mathcal{V},\mathcal{E})$ where $\mathbf{A}_{\mathcal{G}}$ is adjacency matrix whose entry $a_{uv}$ is given by
$$a_{uv}=\left\{\begin{aligned}&1, \quad\text{vertex $u$ is connected to $v$ by an edge}\\
&0,\quad\text{otherwise}.
\end{aligned}\right.
$$
and $\mathbf{D}_{\mathcal{G}}$ is the diagonal matrix that is defined as below: the $u$-th diagonal entry is $k_{u}$, while all non-diagonal entries are zero, i.e., $\mathbf{D}_{\mathcal{G}}=\text{diag}[k_{1},k_{2},\dots,k_{|\mathcal{V}|}]$. Based on this, the hitting time from vertex $u$ to $v$ denoted by $\mathcal{H}_{u\rightarrow v}$, as one significant structural parameter related to random walks, is written as

\begin{equation}\label{eqa:MF-3-0-1}
\\mathcal{H}_{u\rightarrow v}=\left.\dfrac{d}{dz}\mathcal{P}_{u\rightarrow v}(z)\right|_{z=1},
\end{equation}
where $\mathfrak{P}_{u\rightarrow v}(z)$ is the probability generating function that is referred to as

$$\mathfrak{P}_{u\rightarrow v}(z)=\sum_{t=0}^{\infty}P_{u\rightarrow v}(t)z^{t}.$$
Here, we have made use of the probability $P_{u\rightarrow v}(t)$ that a walker starting out from source vertex $u$ first reaches to his destination vertex $v$ after $t$ steps. For network $\mathcal{G}(\mathcal{V},\mathcal{E})$ as a whole, mean hitting time, denoted by $\langle\mathcal{H}_{\mathcal{G}}\rangle$, is viewed as

\begin{equation}\label{eqa:MF-3-0-2}
\langle\mathcal{H}_{\mathcal{G}}\rangle=\frac{1}{|\mathcal{V}|(|\mathcal{V}|-1)}\sum_{u,v}\mathcal{H}_{u\rightarrow v}.
\end{equation}
Note also that the quantity $\mathcal{H}_{u\rightarrow v}$ is not necessarily equal to $\mathcal{H}_{v\rightarrow u}$ in general.

\textbf{Lemma 1 \cite{Zhu-1996}} Given a graph $\mathcal{G}(\mathcal{V},\mathcal{E})$, the exact solution of mean hitting time $\langle\mathcal{H}_{\mathcal{G}}\rangle$ is given by
\begin{equation}\label{eqa:MF-3-0-2-1}
\langle\mathcal{H}_{\mathcal{G}}\rangle=\frac{2|\mathcal{E}|}{|\mathcal{V}|-1}\sum_{i=2}^{|\mathcal{V}|}\frac{1}{\lambda_{i}}
\end{equation}
in which $\lambda_{i}$ is all the non-zero eigenvalues of Laplacian matrix $\mathbf{L}_{\mathcal{G}}(=\mathbf{D}_{\mathcal{G}}-\mathbf{A}_{\mathcal{G}})$ associated with network $\mathcal{G}(\mathcal{V},\mathcal{E})$. 

Suppose that graph in question is a tree $\mathcal{T}$, the summation over quantities $\mathcal{H}_{u\rightarrow v}$ and $\mathcal{H}_{v\rightarrow u}$ is closely connected to their own distance $d_{uv}$, which is shown in the following term

\begin{equation}\label{eqa:MF-3-0-3}
\mathcal{H}_{u\rightarrow v}+\mathcal{H}_{v\rightarrow u}=2d_{uv}(|\mathcal{T}|-1).
\end{equation}
What's more, a more general version associated with Eq.(\ref{eqa:MF-3-0-3}) has been reported in \cite{Chandra-1996,Georgakopoulos-2017}. In fact, such a summation on left-hand side of Eq.(\ref{eqa:MF-3-0-3}) is the so-called commute time between vertices $u$ and $v$. And then, there is a close relation between mean hitting time $\langle\mathcal{H}_{\mathcal{T}}\rangle$ and Wiener index $\mathcal{W}_{\mathcal{T}}$ in the following form.

\textbf{Lemma 2 \cite{MF-2020}} Given a tree $\mathcal{T}$, the closed-form solution of mean hitting time $\langle\mathcal{H}_{\mathcal{T}}\rangle$ is given by
\begin{equation}\label{eqa:MF-3-0-4}
\langle\mathcal{H}_{\mathcal{T}}\rangle=2\mathcal{W}_{\mathcal{T}}/|\mathcal{T}|
\end{equation}
where $\mathcal{W}_{\mathcal{T}}$ is Wiener index of tree $\mathcal{T}$, and is by definition expressed as

$$\mathcal{W}_{\mathcal{T}}=\frac{1}{2}\sum_{u,v\in \mathcal{T}}d_{uv}.$$

Clearly, determining the exact solution to mean hitting time on tree network is equivalently transformed into the problem of calculating the corresponding Wiener index. This enables calculation of mean hitting time on tree when it is more convenient to derive the latter, and vice verse. As will see, the following derivation of mean hitting time on our tree network $\mathcal{T}_{t}$ is carried out based on Eq.(\ref{eqa:MF-3-0-4}).

\subsection{Mean hitting time}

In this section, we will precisely determine the closed-form solution to mean hitting time on tree networks $\mathcal{T}_{t}$ and $\mathcal{T}^{\star}_{t}$. It should be mentioned that a rigorous derivation of mean hitting time $\langle\mathcal{H}_{t}\rangle$ for random walks on tree $\mathcal{T}_{t}$ is deferred to show in Appendix. Furthermore, we immediately give an approximate solution to $\langle\mathcal{H}^{\star}_{t}\rangle$ on tree $\mathcal{T}^{\star}_{t}$ due to a similar discussion to that for quantity $\langle\mathcal{H}_{t}\rangle$.

\textbf{Theorem 5} The closed-from solution of mean hitting time $\langle\mathcal{H}_{t}\rangle$ on tree network $\mathcal{T}_{t}$ is given by

\begin{equation}\label{eqa:MF-3-1-8-6}
\langle\mathcal{H}_{t}\rangle=\frac{2}{5^{t}(|\mathcal{T}_{0}|-1)+1}\left\{75^{t}\mathcal{W}_{0}-4\times5^{t}\left[\frac{15^{t}-5^{t}}{10}(|\mathcal{T}_{0}|-1)^{2}+\frac{15^{t}-1}{14}(|\mathcal{T}_{0}|-1)\right]\right\}
\end{equation}
in which we have made use of Eq.(\ref{eqa:MF-1-1}). A rigorous proof is deferred to show in Appendix.

As a byproduct, the average shortest path length $\langle\mathcal{W}_{t}\rangle$ of scale-free tree network $\mathcal{T}_{t}$ is by definition derived, which is given by

\begin{equation}\label{eqa:MF-A-1}
\begin{aligned}\langle\mathcal{W}_{t}\rangle&=\frac{\mathcal{W}_{t}}{|\mathcal{T}_{t}|(|\mathcal{T}_{t}|-1)/2}\\
&=\frac{2}{5^{t}(|\mathcal{T}_{0}|-1)[5^{t}(|\mathcal{T}_{0}|-1)+1]}\left\{75^{t}\mathcal{W}_{0}-4\times5^{t}\left[\frac{15^{t}-5^{t}}{10}(|\mathcal{T}_{0}|-1)^{2}+\frac{15^{t}-1}{14}(|\mathcal{T}_{0}|-1)\right]\right\}.
\end{aligned}
\end{equation}
In the thermodynamic limit, the scaling of quantity $\langle\mathcal{W}_{t}\rangle$ obeys

\begin{equation}\label{eqa:MF-A-2}
\langle\mathcal{W}_{t}\rangle\sim3^{t}=O(|\mathcal{T}_{t}|^{\zeta}), \quad \text{and},\quad \zeta=\ln3/\ln5.
\end{equation}

Clearly, average shortest path length $\langle\mathcal{W}_{t}\rangle$ and diameter $D_{t}$ in tree network $\mathcal{T}_{t}$ have the same order of magnitude. In other words, we show using result in Eq.(\ref{eqa:MF-A-2}) that tree network $\mathcal{T}_{t}$ indeed has no small-world property.

It is worth mentioning that different from most prior work focusing particularly on tree networks whose seed is a single edge, formula in Eq.(\ref{eqa:MF-3-1-8-6}) holds in the sense that an arbitrary tree is selected as the seed for creating tree network $\mathcal{T}_{t}$. Additionally, as an illustrative example, we plot numerical values for $\langle\mathcal{H}_{t}\rangle$ at different time steps in Fig.3 where four simple trees, i.e., an edge,  a path of length $3$, a path of length $5$ and a star with $6$ vertices, are selected to serve as seed, respectively. The results suggest that empirical simulations are in perfect agreement with the theoretical analysis.

\begin{figure*}
\centering
\subfigure[An edge as seed]{
\begin{minipage}[t]{0.5\linewidth}
\centering
\includegraphics[width=7.5cm]{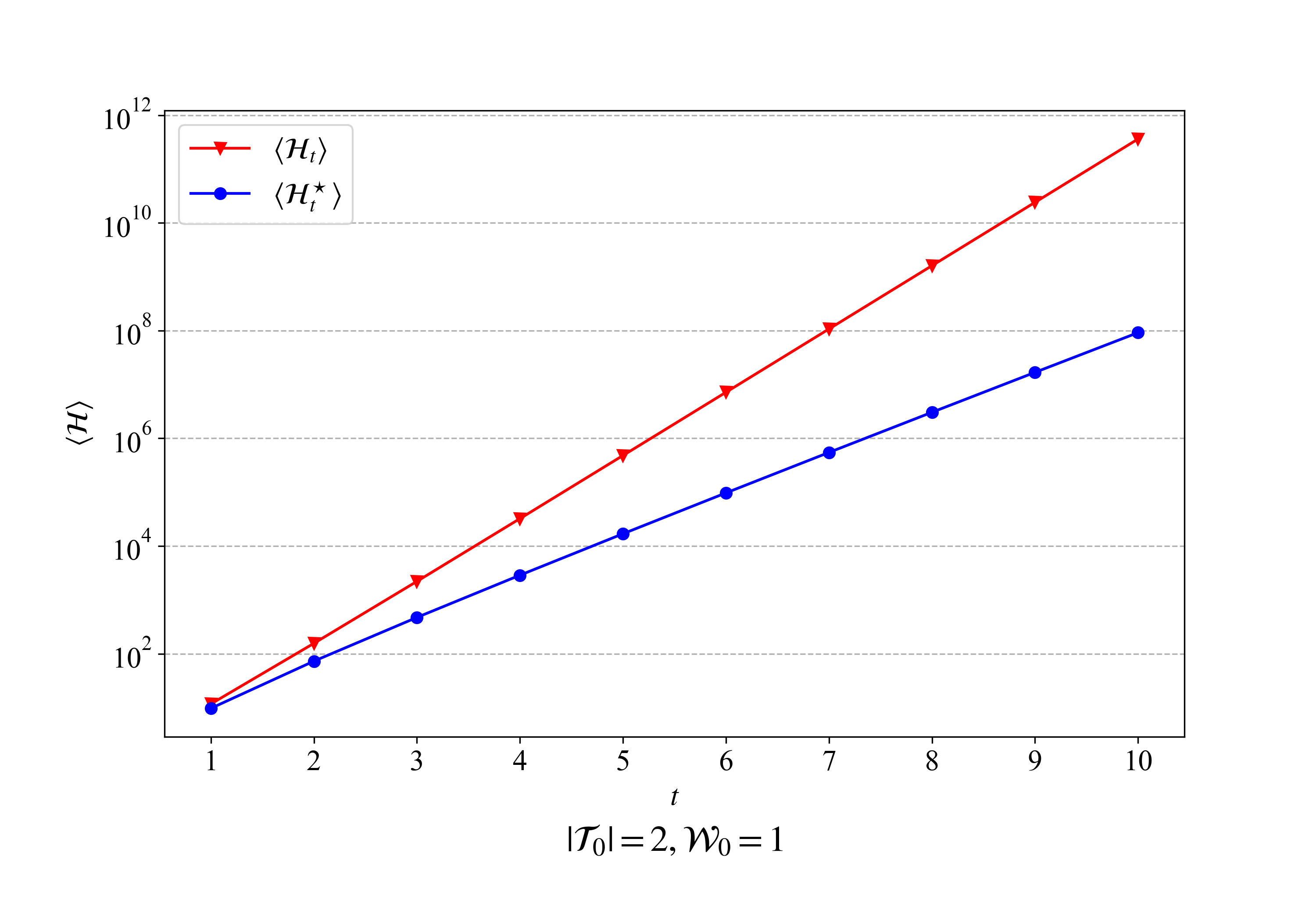}
%\caption{fig1}
\end{minipage}
}%
\subfigure[A path of length $3$ as seed]{
\begin{minipage}[t]{0.5\linewidth}
\centering
\includegraphics[width=7.5cm]{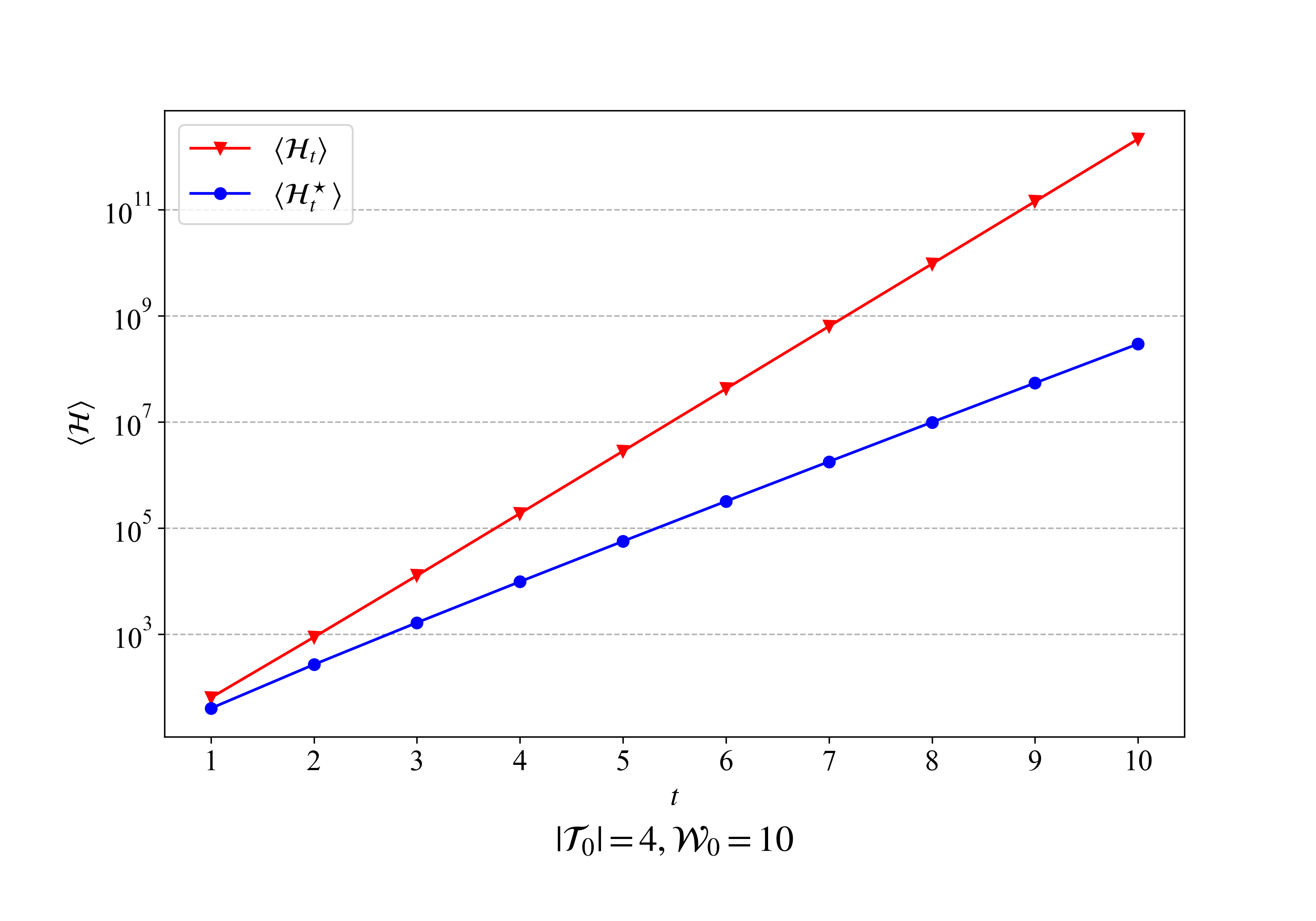}
%\caption{fig2}
\end{minipage}
}%
\\
\subfigure[A path of length $5$ as seed]{
\begin{minipage}[t]{0.5\linewidth}
\centering
\includegraphics[width=7.5cm]{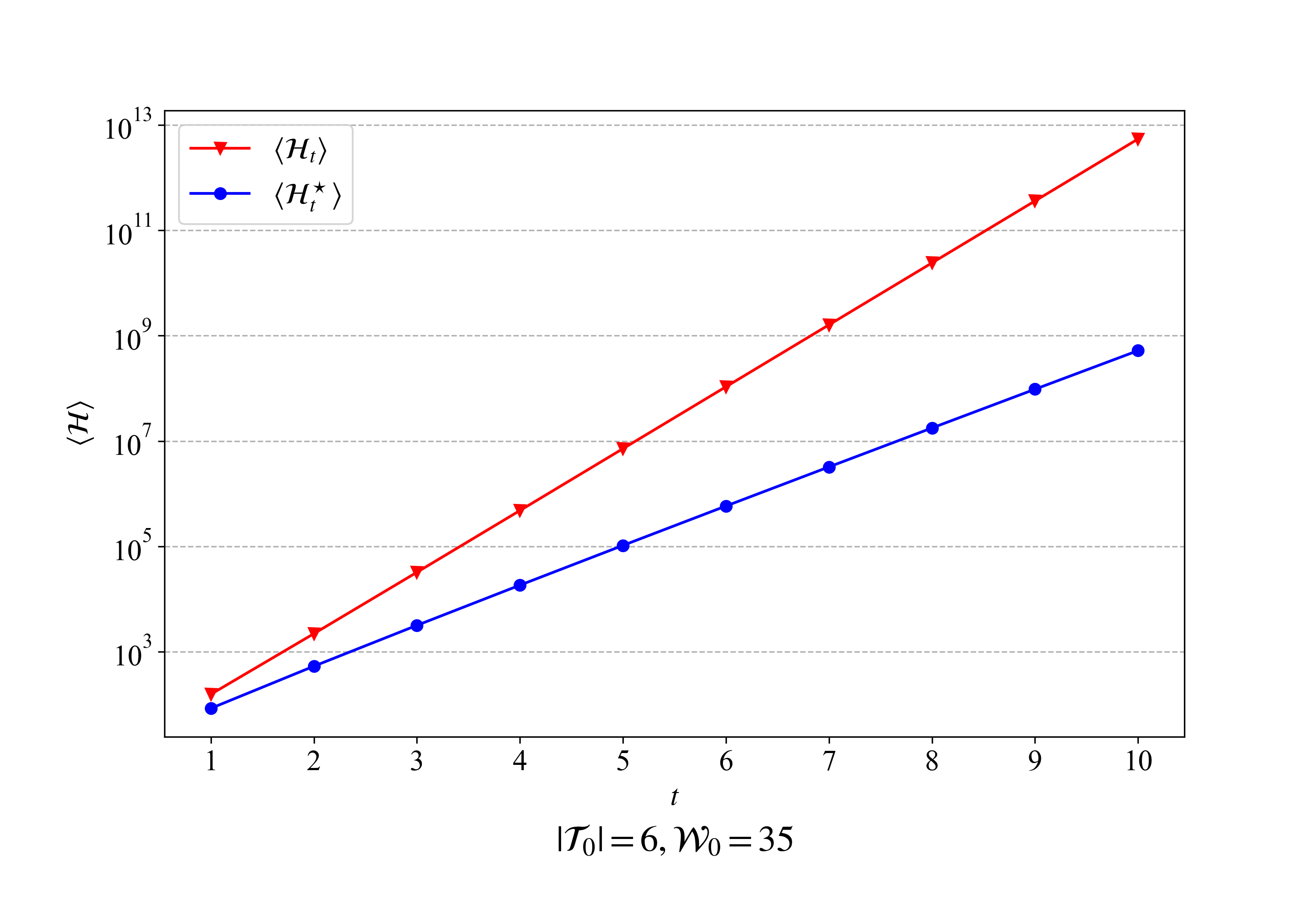}
%\caption{fig2}
\end{minipage}
}%
\subfigure[A star with $6$ vertices as seed]{
\begin{minipage}[t]{0.5\linewidth}
\centering
\includegraphics[width=7.5cm]{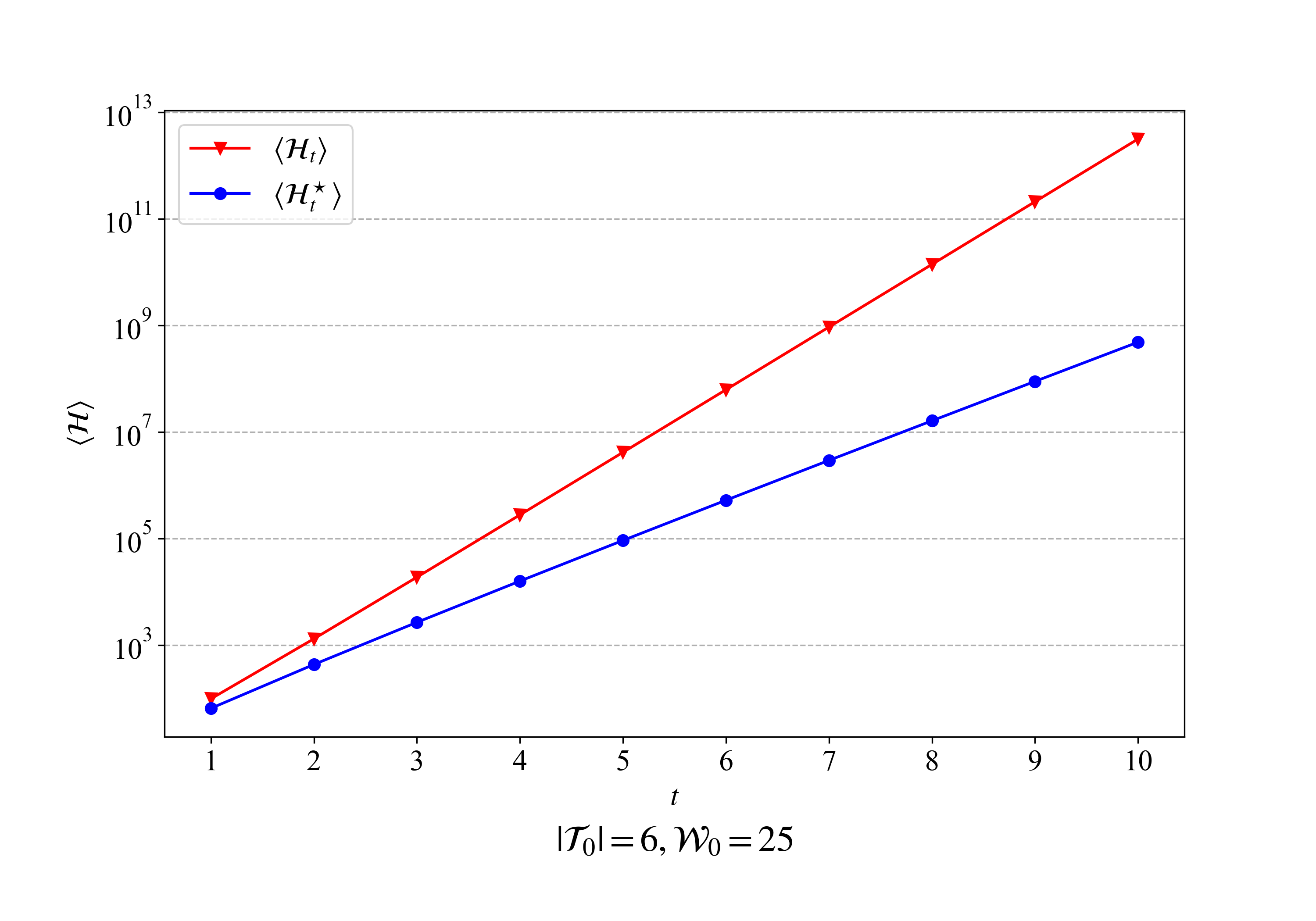}
%\caption{fig3}
\end{minipage}
}%
\caption{(Color online)  The diagram of mean hitting time $\langle\mathcal{H}_{t}\rangle$ and $\langle\mathcal{H}_{t}^{\star}\rangle$. Note that we select four simple trees, i.e., an edge, a path of length $3$, a path of length $5$ and a star with $6$ vertices, as the seed for creating tree networks. In all the panels, the straight line indicates analytical solutions for quantities $\langle\mathcal{H}_{t}\rangle$ and $\langle\mathcal{H}_{t}^{\star}\rangle$. Accordingly, triangle and circle represent the corresponding numerical values produced by computer simulation. Clearly, theoretical analysis is strongly consistent with experimental simulation. At the same time, quantities $\langle\mathcal{H}_{t}\rangle$ and $\langle\mathcal{H}_{t}^{\star}\rangle$ have different scaling phenomena in the limit of large graph size. Also, $\langle\mathcal{H}_{t}\rangle$  is strictly larger than $\langle\mathcal{H}_{t}^{\star}\rangle$. }
\end{figure*}

In the thermodynamic limit, we can find the next expression

\begin{equation}\label{eqa:MF-3-1-8-7}
\langle\mathcal{H}_{t}\rangle\sim|\mathcal{T}_{t}|^{\chi}, \quad \chi=1+\frac{\ln3}{\ln5}.
\end{equation}
In the past, tree models whose mean hitting time exhibits a power-law function about vertex number are rarely encountered in the study of scale-free tree networks. Although some fractal trees, for instance, Vicsek fractal \cite{Vicsek-1983} and T-graph \cite{Redner-2001}, show a similar phenomenon when studying mean hitting time on their own underlying structures, it is well known that they have no scale-free feature. In a nutshell, our tree network $\mathcal{T}_{t}$, as one of important members in tree family, helps one to in depth probe rich underlying structures of tree models. At the same time, our analysis also provides some instruments for building some tree networks of great interest.

As mentioned previously, mean hitting time $\langle\mathcal{H}^{\star}_{t}\rangle$ for random walks on scale-free tree network $\mathcal{T}^{\star}_{t}$ is exactly derived in a similar manner as well. In fact, almost all calculations necessary for deriving quantity  $\langle\mathcal{H}^{\star}_{t}\rangle$ have been reported in the process of computing $\langle\mathcal{H}_{t}\rangle$ and one just needs to make some slight modifications before reaching his goal. Therefore, we omit the corresponding detailed derivations and only provide the final results as below

\begin{equation}\label{eqa:MF-3-1-8-8-1}
\langle\mathcal{H}^{\star}_{t}\rangle=\frac{2}{5^{t}(|\mathcal{T}_{0}|-1)+1}\left[25^{t}\mathcal{W}_{0}+10t\times5^{2(t-1)}(|\mathcal{T}_{0}|-1)^{2}-\frac{3\times(5^{2t-1}-5^{t-1})}{2}(|\mathcal{T}_{0}|-1)\right].
\end{equation}
Yet, numerical simulations for mean hitting time $\langle\mathcal{H}^{\star}_{t}\rangle$ are plotted in Fig.3. 

Clearly, we have an approximate solution to quantity $\langle\mathcal{H}^{\star}_{t}\rangle$, i.e.,

\begin{equation}\label{eqa:MF-3-1-8-8}
\langle\mathcal{H}^{\star}_{t}\rangle\sim|\mathcal{T}^{\star}_{t}|^{\chi^{\star}}\ln|\mathcal{T}^{\star}_{t}|, \quad \chi^{\star}=1
\end{equation}
where we have used Eq.(\ref{eqa:MF-1-2}). It should be mentioned that quantity $\langle\mathcal{H}^{\star}_{t}\rangle$ is not studied in \cite{Jung-2002}. From Eqs.(\ref{eqa:MF-1-1}),(\ref{eqa:MF-1-2}),(\ref{eqa:MF-3-1-8-7}) and (\ref{eqa:MF-3-1-8-8}), we see that while trees $\mathcal{T}_{t}$ and $\mathcal{T}^{\star}_{t}$ have the same vertex number, their own underlying structures play completely different roles on studying random walks.

\subsection{Discussions}

Below, we continue to discuss about trees $\mathcal{T}_{t}$ and $\mathcal{T}^{\star}_{t}$ in more detail for the purpose of revealing similarities and differences between these two tree networks. In addition, some other related topics are also considered. That is to say, this section aims at extending some results in the previous sections.

First of all, let us rethink of diameter and mean hitting time on two tree networks $\mathcal{T}_{t}$ and $\mathcal{T}^{\star}_{t}$. In the large graph size limit, there is an apparent connection, namely,

$$\langle\mathcal{H}_{t}\rangle/|\mathcal{T}_{t}|\sim D_{t}, \quad \langle\mathcal{H}^{\star}_{t}\rangle/|\mathcal{T}^{\star}_{t}|\sim D^{\star}_{t}.$$
This tells us that for a tree network $\mathcal{T}$, diameter can be used as an index for approximately measuring mean hitting time for random walks. One of important reasons for this is the elegant relationship between Wiener index $\mathcal{W}_{\mathcal{T}}$ and mean hitting time $\langle\mathcal{H}_{\mathcal{T}}\rangle$ mentioned in Eq.(\ref{eqa:MF-3-0-4}). Roughly speaking, if one can determine a precise enough solution of diameter of tree network, then it is not necessary to calculate the corresponding mean hitting time when he only wants an approximate description about quantity $\langle\mathcal{H}_{\mathcal{T}}\rangle$.

On the other hand, from the theory point of view, it is interesting yet challenging to capture exact solution of mean hitting time on various kinds of trees. It is mainly because in a rigorous calculational process, one may develop some effective techniques in order to address some tricky tasks and also uncover connections between some topological parameters behind tree networks. For instance, the term $\sum_{v\in\mathcal{T}_{t-1}}[k_{v}(t-1)]^{2}$ in Eq.(\ref{subeq:MF-3-1-3-3}) is not derived straightforwardly. To handle this issue, we establish convenient transformation methods from unknowns to some quantities known to us. As a result, some simple and elegant relationships between distinct quantities are built, such as, formula in Eq.(\ref{eqa:MF-3-1-8-5}). This further enables understanding of topological structure of tree networks. In a word, it is still of great interest to derive closed-form solutions of some fundamental and important structural parameters on a variety of tree networks in a mathematically rigorous manner.

As exhibited by Algorithm-I in Section 2, tree network $\mathcal{T}_{t}$ is generated via two different kinds of mechanisms for adding new vertices. However, only one of both mechanisms is used to create tree network $\mathcal{T}^{\star}_{t}$. Our analysis shows that doing so leads to two different kinds of trees. Although they have scale-free feature, the power-law exponents fall into two distinct regions. As pointed out in our prior work \cite{Ma-2020-T}, if one only inserts new vertices on existing edge of tree, the mean hitting time for the end tree can asymptotically reach to the theoretical upper bound in the large graph size limit. On the contrary, connecting new vertices to each vertex of tree directly at each time step leads to a tree whose mean hitting time gradually tends to the opposite direction, i.e., to the theoretical lower bound. Such an example is tree $\mathcal{T}^{\star}_{t}$. Now, it is clear to see that the scaling of mean hitting time on tree network $\mathcal{T}_{t}$ is within the middle region between theoretical lower bound and theoretical upper bound. While some fractal trees \cite{Vicsek-1983}-\cite{Agliari-2008} possess such a characteristic, they do not have scale-free feature. On balance, the light shed by constructing our tree $\mathcal{T}_{t}$ is helpful to generate other interesting trees in the future.

\begin{figure*}
\centering
\subfigure[An edge as seed]{
\begin{minipage}[t]{0.5\linewidth}
\centering
\includegraphics[width=8cm]{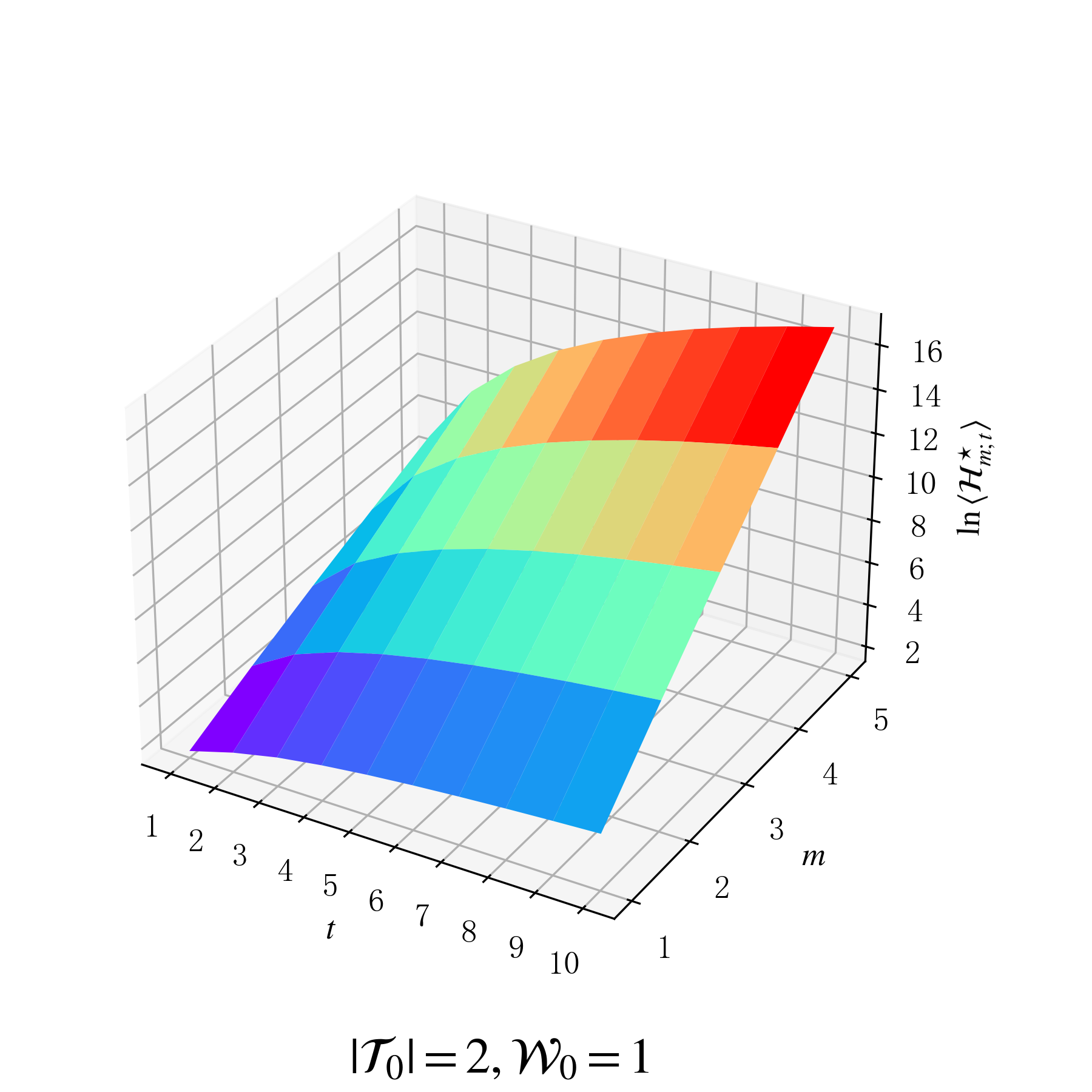}
%\caption{fig1}
\end{minipage}
}%
\subfigure[A path of length $3$ as seed]{
\begin{minipage}[t]{0.5\linewidth}
\centering
\includegraphics[width=8cm]{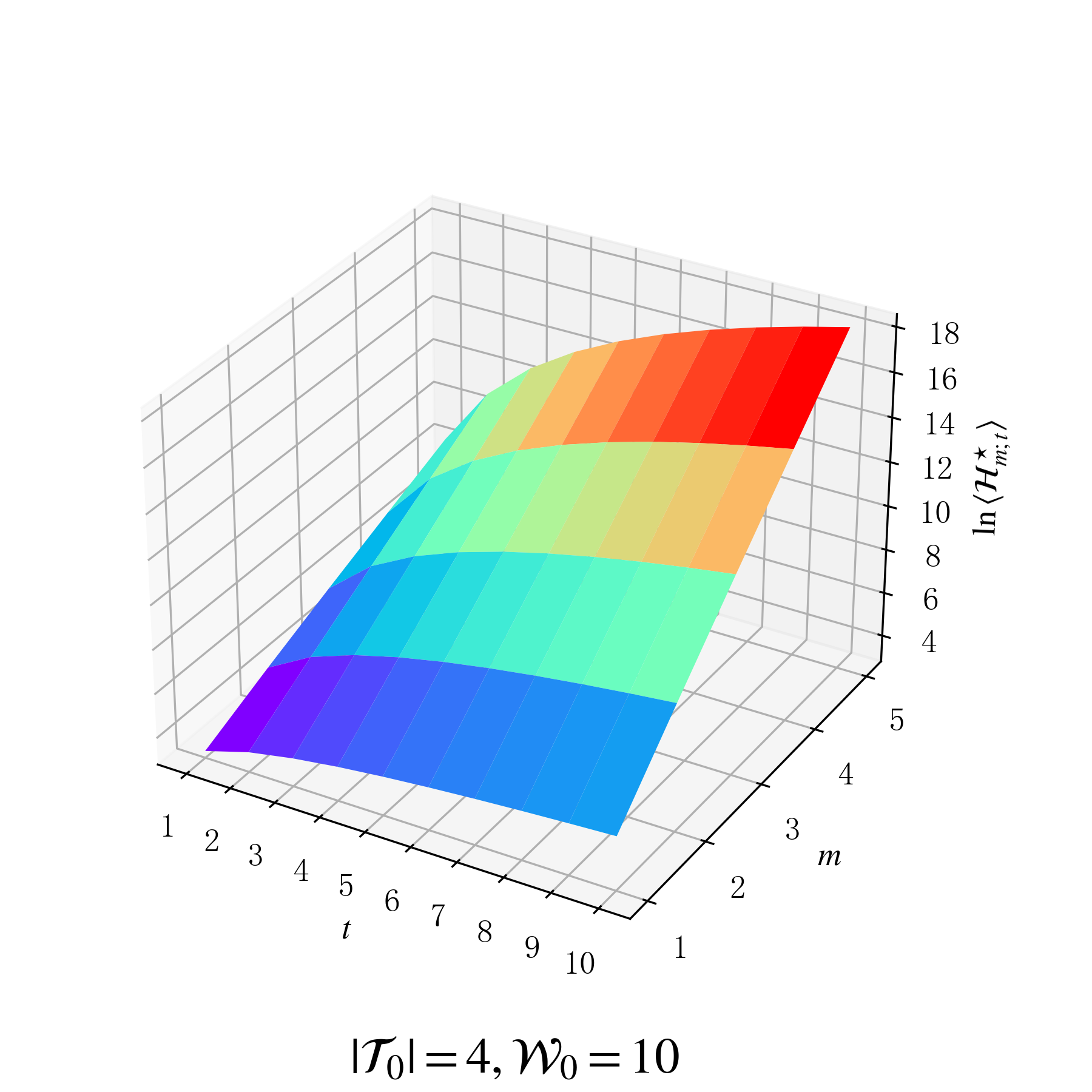}
%\caption{fig2}
\end{minipage}
}%
\\
\subfigure[A path of length $5$ as seed]{
\begin{minipage}[t]{0.5\linewidth}
\centering
\includegraphics[width=8cm]{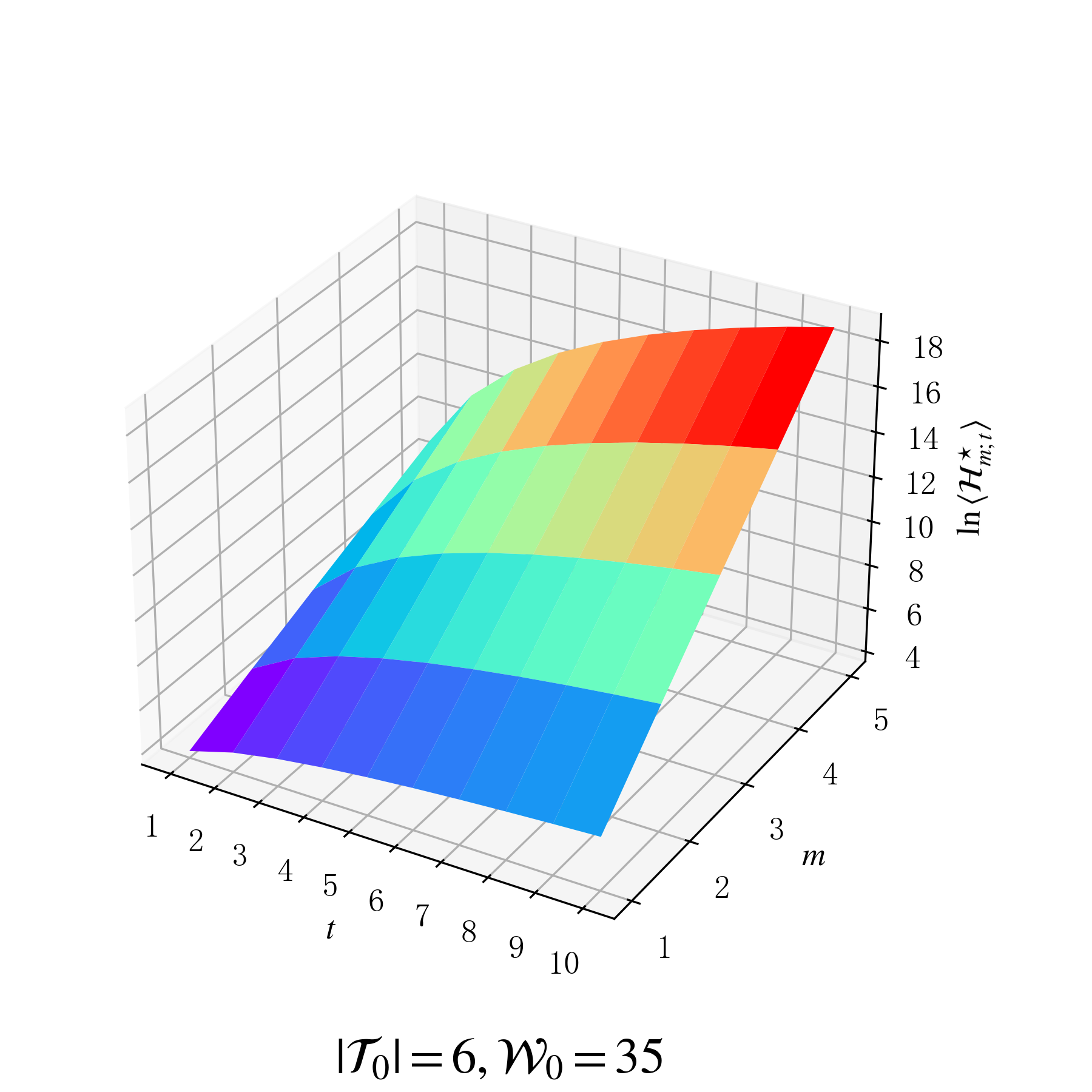}
%\caption{fig2}
\end{minipage}
}%
\subfigure[A star with $5$ vertices as seed]{
\begin{minipage}[t]{0.5\linewidth}
\centering
\includegraphics[width=8cm]{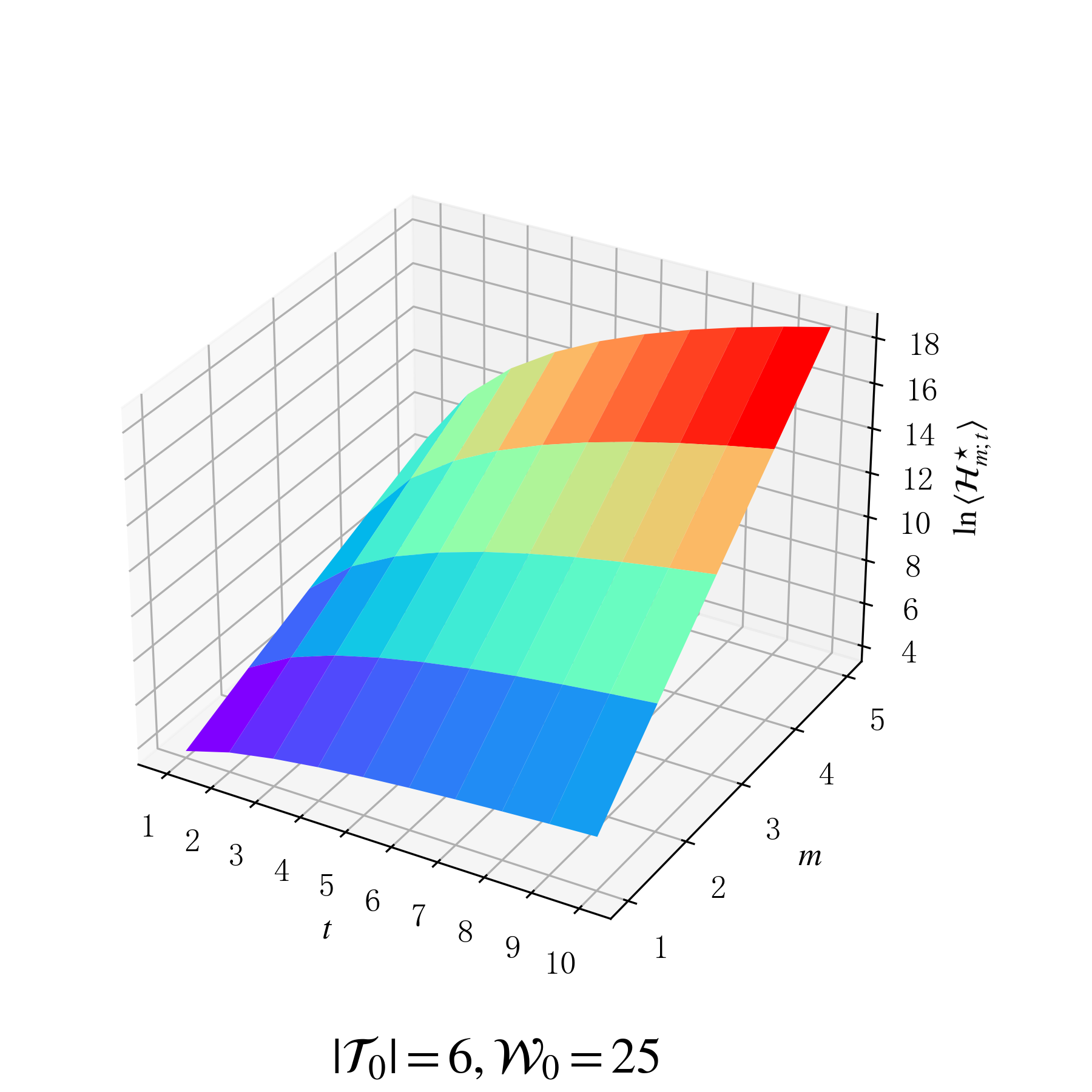}
%\caption{fig3}
\end{minipage}
}%
\caption{(Color online)  The diagram of mean hitting time $\langle\mathcal{H}^{\star}_{m;t}\rangle$. Note that we still select four simple trees, i.e., an edge, a path of length $3$, a path of length $5$ and a star with $6$ vertices, as the seed for creating tree networks. In all the panels, the straight line indicates analytical solution. Accordingly, square represents the corresponding numerical value produced by computer simulation. Clearly, experimental simulations are strongly consistent with theoretical analysis. }
\end{figure*}

At last, let us generalize the procedure introduced in Algorithm-I in order to create other scale-free tree networks with both power-law exponent larger than $3$ and larger diameter. Concretely speaking, the generalization also consists of two ingredients at each time step. The first is to insert $n$ new vertices on each existing edge $uv$. The other is to connect $mk_{v}$ new vertices as leaves to each pre-existing vertex $v$. After $t$ time steps, the resultant tree network is denoted by $\mathcal{T}_{m,n;t}$ for convenience. Next, we can employ a similar method as used to prove theorem 4.1 to verify that tree network $\mathcal{T}_{m,n;t}$ follows power-law degree distribution with exponent $\ln(n+2m+1)/\ln(m+1)$. It is clear to see that if parameters $n$ and $m$ satisfy an inequality $(n+2m+1)>(m+1)^{3}$, the resulting scale-free tree network $\mathcal{T}_{m,n;t}$ has power-law exponent strictly larger than $3$. At the same time, it is not hard to obtain that diameter of scale-free tree network $\mathcal{T}_{m,n;t}$ is subject to $O((n+1)^{t})$. These such scale-free tree networks are what we want. In addition, we may also derive the solution of mean hitting time $\langle\mathcal{H}_{m,n;t}\rangle$ for random walks on tree networks $\mathcal{T}_{m,n;t}$ in a manner established in proof of theorem 5. For the sake of space, we only give the asymptotical solution of mean hitting time, i.e., $\langle\mathcal{H}_{m,n;t}\rangle\sim|\mathcal{T}_{m,n;t}|^{1+\frac{\ln(n+1)}{\ln(n+2m+1)}}$. If $n\gg m$, we have $\langle\mathcal{H}^{\star}_{m;t}\rangle\sim O\left(|\mathcal{T}_{m,n;t}|^{2}\right)$. In this setting, scale-free tree network $\mathcal{T}_{m,n;t}$ nearly achieves theoretical upper bound on mean hitting time \cite{Hwang-2012}.

Along the similar research line as above, we may generate scale-free tree network with small diameter as well. Specifically speaking, we only need to connect $mk_{v}$ new vertices as leaves to each pre-existing vertex $v$ at each time step. The $t$-th generation is denoted by $\mathcal{T}^{\star}_{m;t}$. Accordingly, some fundamental structural parameters on tree network $\mathcal{T}^{\star}_{m;t}$ are analyzed in similar manner as used in the previous sections. In the following, we just provide the solution of mean hitting time 

$$
\langle\mathcal{H}^{\star}_{m;t}\rangle=\frac{2\left\{[f_{1}(m)]^{t}\mathcal{W}_{0}+f_{2}(m)\sum_{i=0}^{t-1}[f_{1}(m)]^{i}|\mathcal{T}^{\star}_{m;t-1-i}|^{2}-f_{3}(m)\sum_{i=0}^{t-1}[f_{1}(m)]^{i}|\mathcal{T}^{\star}_{m;t-1-i}|+f_{4}(m)\sum_{i=0}^{t-1}[f_{1}(m)]^{i}\right\}}{|\mathcal{T}^{\star}_{m;t}|},
$$
in which $f_{1}(m)=(2m+1)^{2}$, $f_{2}(m)=m(2m+1)$, $f_{3}(m)=m(5m+3)$, $f_{4}(m)=m(3m+2)$ and $|\mathcal{T}^{\star}_{m;t}|=(|\mathcal{T}_{0}|-1)(2m+1)^{t}+1$. To make our statement more concrete, numerical analysis of quantity $\langle\mathcal{H}^{\star}_{m;t}\rangle$ is plotted in Fig.4. It is obvious that theoretical analysis is in strong agreement with experimental simulations.

\section{Related work}

Scale-free networks have attracted more attention since the famous BA-model attributed to Barab\'{a}si and Albert \cite{Barabasi-1999}-\cite{Telesford-2011}. As a result, a great number of scale-free network models have been proposed, including static scale-free network model \cite{Catanzaro-2005}, scale-free networks without growth \cite{Xie-2008}, vertex-edge-growth scale-free network models \cite{Ma-2018}, Apollonian networks \cite{Andrade-2005}, pseudofractal scale-free web \cite{Dorogovtsev-2002}, hierarchical networks \cite{Ravasz-2003} and flower graphs \cite{Diggans-2020}. Accordingly, these previously published scale-free networks have analyzed in depth by studying some fundamental structural parameters such as degree distribution, diameter etc. As shown in \cite{Cohen-2003}, Cohen \emph{et al} have shown that scale-free networked models are small and, sometimes, ultra-small by determining relationship between diameter and vertex number. On the other hand, the prior work shows that there are theoretical scale-free models possessing a larger diameter \cite{Zhang-2007}-\cite{Wang-2020-chaos}. This suggests that there is room to completely understand underlying structure of scale-free networks. In another word, it is of great interest to establish new scale-free networks in the current study. 

Tree network, as the simplest connected network, plays a key role in the field of scale-free networks \cite{Athreya-2017}-\cite{Ma-2020-T}. Among of them, the well-known is BA-tree induced by famous BA-model \cite{Barabasi-1999}. In \cite{Souza-2007}, Souza \emph{et al} proposed a family of random recursive trees and found that greedy choice broadens the degree distribution according to the power of choice. Additionlly, some deterministic trees with intriguing structural properties have also been established, such as, geometric fractal growth tree model \cite{Jung-2002}, deterministic uniform recursive tree \cite{Lu-2010}, growing treelike networks \cite{Zhang-2010}, Fibonacci treelike models \cite{Ma-2020-1}, Vicsek fractal \cite{Vicsek-1983}-\cite{Ma-epl-2021} and T-graph \cite{Redner-2001}-\cite{Ma-2020-T}. The previous results have shown that scale-free feature and larger diameter can not be found on these pre-existing tree networks simultaneously. This leads to an urgent need to create tree networks where the above-mentioned two properties coexist.

After building up networks, it is sufficient to investigate the associated underlying structure by determining some fundamental parameters in order to understand them in more detail. Sometimes, some parameters are introduced by considering dynamics such as random walks taking place on networks \cite{Dyer-2020}-\cite{Mokhtar-2013} and \cite{Hwang-2012,Perkins-2014}. For example, mean hitting time (also called mean first-passage time in some published papers) for random walks on networks including scale-free models has received in the past \cite{Vicsek-1983,Zhang-2010}. Zhang \emph{et al} \cite{Zhang-2010} derived mean first-passage time for random walks on growing treelike networks using spectral technique. The similar method has also been used to obtain hitting time on other networks \cite{Dyer-2020,Noh-2004}. More generally, it is challenging to determine exact solution to some structural parameters including mean hitting time on a growth tree network with thousands of vertices. Therefore, it is necessary to propose proper technique when considering some specific networks.

\section{Conclusion}

In summary, we put forward a class of growth tree networks $\mathcal{T}_{t}$ and discuss some fundamental structural parameters. Model $\mathcal{T}_{t}$ turns out to follow power-law degree distribution with exponent larger than $3$, which is seldom reported in the realm of scale-free tree networks. Next, small-world property can not be observed on tree $\mathcal{T}_{t}$ because it has a large diameter, which is completely different from most of published scale-free trees. We also show that tree network $\mathcal{T}_{t}$ possesses disassortative structure. Additionally, in order to better understand the underlying structure of tree $\mathcal{T}_{t}$, we study random walks and then derive an exact formula for mean hitting time. The results show that as opposed to most previously proposed scale-free trees, the mean hitting time of our tree network $\mathcal{T}_{t}$ obeys a power-law function as vertex number of tree $\mathcal{T}_{t}$. Note that another scale-free network $\mathcal{T}^{\star}_{t}$ is in depth analyzed as well. As mentioned in the text, tree network $\mathcal{T}_{t}$ shows significantly distinct topological structure with its analog $\mathcal{T}^{\star}_{t}$ in terms of several structural indices while they have a vertex number in common. Last but not the least, we firmly believe that tree network $\mathcal{T}_{t}$ has many other intriguing features that are worthy of attention, which is left as research direction in our future work.

\section*{Acknowledgments}

We would like to thank anonymous reviewers for their valuable comments on our paper, which have considerably improved the presentation of this paper. The research was supported by the Fundamental Research Funds for the Central Universities No. G2023KY05105 and the National Key Research and Development Plan under grant 2020YFB1805400.

\section*{Appendix}
\setcounter{equation}{0}                                                                                                                                                                                                                                                                               
\renewcommand\theequation{A.\arabic{equation}}

Here is a proof of theorem 5. 

\emph{Proof } Before beginning with our calculations, let us recall the concrete procedure for developing scale-free tree network $\mathcal{T}_{t}$. At time step $t$, tree $\mathcal{T}_{t}$ consists of three classes of vertex sets, i.e., $\mathcal{T}_{t-1}$, $A(t)=\bigcup_{v\in\mathcal{T}_{t-1}}A_{v}(t)$ and $B(t)=\bigcup_{v\in\mathcal{T}_{t-1}}B_{v}(t)$. It is natural to see that by definition, Wiener index $\mathcal{W}_{t}$ of tree $\mathcal{T}_{t}$ is obtained if we can determine distances of all possible pairs of vertices from the three sets. Roughly speaking, there are nine distinct combinatorial manners according to different classifications of vertices of tree $\mathcal{T}_{t}$. For brevity and our purpose, given a vertex $v$ of tree $\mathcal{T}_{t-1}$, each vertex in set $A_{v}(t)$ is labeled as $v^{1}_{i}$ and, similarly, each vertex in set $B_{v}(t)$ is labeled as $v^{2}_{i}$. Now, we are ready to estimate quantity $\langle\mathcal{H}_{t}\rangle$ on tree $\mathcal{T}_{t}$.

\emph{Stage 1.} Due to the fact that three are only two new vertices inserted on each edge $uv$ of tree $\mathcal{T}_{t-1}$, we can without difficulty obtain

\begin{equation}\label{eqa:MF-3-1-1}
\begin{aligned}\mathcal{W}_{t}(1)&:=\frac{1}{2}\sum_{u,v\in \mathcal{T}_{t-1}}d^{\dagger}_{uv}\\
&=3\mathcal{W}_{t-1}.
\end{aligned}
\end{equation}
Here, $d^{\dagger}_{uv}$ represents distance between vertices $u$ and $v$ when considering tree $\mathcal{T}_{t}$. It is worth noticing that, hereafter, we make use of the superscribe $\dagger$ to indicate distance between vertex pair on tree $\mathcal{T}_{t}$.

\emph{Stage 2.} For each vertex $v$ in tree $\mathcal{T}_{t-1}$, there are in practice $2k_{v}(t-1)$ vertices added into tree $\mathcal{T}_{t}$ where $k_{v}(t-1)$ is the degree of vertex $v$ in tree $\mathcal{T}_{t-1}$. Half of them are grouped into set $A_{v}(t)$. The others are in set $B_{v}(t)$. Using notations defined above, it is easy to see the following two equations

\begin{subequations}
\label{eq:whole}
\begin{eqnarray}
\begin{aligned}\mathcal{W}_{t}(v;1)&:=\sum_{v^{1}_{i}\in A_{v}(t)}d^{\dagger}_{vv^{1}_{i}}\\
&=k_{v}(t-1),
\end{aligned}\label{subeq:MF-3-1-2-1}
\end{eqnarray}
\begin{equation}
\begin{aligned}\mathcal{W}_{t}(v;2)&:=\sum_{v^{2}_{i}\in B_{v}(t)}d^{\dagger}_{vv^{2}_{i}}\\
&=k_{v}(t-1).
\end{aligned}\label{subeq:MF-3-1-2-2}
\end{equation}
\end{subequations}
In a word, we have

\begin{equation}\label{eqa:MF-3-1-2}
\begin{aligned}\mathcal{W}_{t}(2)&:=\sum_{v\in\mathcal{T}_{t-1}}\sum_{v^{1}_{i}\in A_{v}(t)}d^{\dagger}_{vv^{1}_{i}}+\sum_{v\in\mathcal{T}_{t-1}}\sum_{v^{2}_{i}\in B_{v}(t)}d^{\dagger}_{vv^{2}_{i}}\\
&=\sum_{v\in\mathcal{T}_{t-1}}\left(\mathcal{W}_{t}(v;1)+\mathcal{W}_{t}(v;2)\right)\\
&=2\sum_{v\in\mathcal{T}_{t-1}}k_{v}(t-1)\\
&=4(|\mathcal{T}_{t-1}|-1).
\end{aligned}
\end{equation}

\emph{Stage 3.} Following the above-mentioned stage, let us still concentrate on local structure associated with each vertex $v$ in tree $\mathcal{T}_{t-1}$, i.e., $A_{v}(t)$ and $B_{v}(t)$. And then, it is easy to find

\begin{subequations}
\label{eq:whole}
\begin{eqnarray}
\begin{aligned}\mathcal{W}_{t}(A_{v}(t))&:=\frac{1}{2}\sum_{v^{1}_{i},v^{1}_{j}\in A_{v}(t)}d^{\dagger}_{v^{1}_{i}v^{1}_{j}}\\
&=k_{v}(t-1)(k_{v}(t-1)-1),
\end{aligned}\label{subeq:MF-3-1-3-1}
\end{eqnarray}
\begin{equation}
\begin{aligned}\mathcal{W}_{t}(B_{v}(t))&:=\frac{1}{2}\sum_{v^{2}_{i},v^{2}_{j}\in B_{v}(t)}d^{\dagger}_{v^{2}_{i}v^{2}_{j}}\\
&=k_{v}(t-1)(k_{v}(t-1)-1).
\end{aligned}\label{subeq:MF-3-1-3-2}
\end{equation}
\end{subequations}
Then, we write

\begin{subequations}
\label{eq:whole}
\begin{eqnarray}
\begin{aligned}\mathcal{W}_{t}(3)&:=\sum_{v\in\mathcal{T}_{t-1}}\mathcal{W}_{t}(A_{v}(t))\\
&=\sum_{v\in\mathcal{T}_{t-1}}[k_{v}(t-1)]^{2}-2(|\mathcal{T}_{t-1}|-1),
\end{aligned}\label{subeq:MF-3-1-3-3}
\end{eqnarray}
\begin{equation}
\begin{aligned}\mathcal{W}_{t}(3)&:=\sum_{v\in\mathcal{T}_{t-1}}\mathcal{W}_{t}(B_{v}(t))\\
&=\sum_{v\in\mathcal{T}_{t-1}}[k_{v}(t-1)]^{2}-2(|\mathcal{T}_{t-1}|-1).
\end{aligned}\label{subeq:MF-3-1-3-4}
\end{equation}
\end{subequations}

\emph{Stage 4.} Now, we attempt to study distance $d^{\dagger}_{uv^{1}_{i}}$ between vertex $u$ from tree $\mathcal{T}_{t-1}$ and vertex $v^{1}_{i}$ in set $A_{v}(t)$. Note that vertex $u$ is distinct with $v$. In another word, vertex $v^{1}_{i}$ does not belong to set $A_{u}(t)$. To this end, we need to consider a trivial fact that there is a unique path between an arbitrarily given pair of vertices, say $u$ and $v$, in tree. For convenience, we denote by $\mathcal{P}_{uv}$ the path joining vertex $u$ to $v$. According to underlying structure of tree $\mathcal{T}_{t}$, path $\mathcal{P}_{uv^{1}_{i}}$ is constituted in two manners. The first one is that path $\mathcal{P}_{uv^{1}_{i}}$ consists of two segments $\mathcal{P}_{uv}$ and $\mathcal{P}_{vv^{1}_{i}}$. The number of paths of such type is $(k_{v}(t-1)-1)$ in total. The other is that path $\mathcal{P}_{uv^{1}_{i}}$ is obtained from path $\mathcal{P}_{uv}$ by deleting segment $\mathcal{P}_{v^{1}_{i}v}$. This case is encountered in tree $\mathcal{T}_{t}$ only once. Intuitively, given a pair of vertices $u$ and $v$ in tree $\mathcal{T}_{t-1}$,  it is straightforward to find the following expression

\begin{equation}\label{eqa:MF-3-1-4-1}
\begin{aligned}\mathcal{W}_{t}(u;A_{v}(t))&:=\sum_{v(\neq u)\in\mathcal{T}_{t-1}}\sum_{v^{1}_{i}\in A_{v}(t)}d^{\dagger}_{uv^{1}_{i}}\\
&=\sum_{v(\neq u)\in\mathcal{T}_{t-1}}\left[k_{v}(t-1)d^{\dagger}_{uv}+(k_{v}(t-1)-2)\right].
\end{aligned}
\end{equation}

Summing on both hand sides of Eq.(\ref{eqa:MF-3-1-4-1}) with respect to all vertices $u$ in tree $\mathcal{T}_{t-1}$ yields

\begin{equation}\label{eqa:MF-3-1-4-2}
\begin{aligned}\mathcal{W}_{t}(4)&:=\sum_{u\in\mathcal{T}_{t-1}}\mathcal{W}_{t}(u;A_{v}(t))\\
&=\sum_{u\in\mathcal{T}_{t-1}}\sum_{v(\neq u)\in\mathcal{T}_{t-1}}k_{v}(t-1)d^{\dagger}_{uv}+\sum_{u\in\mathcal{T}_{t-1}}\sum_{v(\neq u)\in\mathcal{T}_{t-1}}(k_{v}(t-1)-2).
\end{aligned}
\end{equation}
However, it is not easy to derive an exact solution for quantity $\mathcal{W}_{t}(4)$ in terms of Eq.(\ref{eqa:MF-3-1-4-2}). To get around this issue, we introduce another calculational technique as below. Specifically, we take an arbitrary path $\mathcal{P}_{uv}$ in tree $\mathcal{T}_{t}$ as an illustrative example. If the length of path $\mathcal{P}_{uv}$ is equal to $3$, we can express $\mathcal{P}_{uv}:=uu^{1}_{i}v^{1}_{j}v$. In this case, it is easy to obtain
\begin{equation}\label{eqa:MF-3-1-4-3}
d^{\dagger}_{uv^{1}_{j}}+d^{\dagger}_{vu^{1}_{i}}=4.
\end{equation}
Due to nature of tree, there are $(|\mathcal{T}_{t-1}|-1)$ paths of this kind in total. Besides that, if path $\mathcal{P}_{uv}$ has length larger than $3$, then, without loss of generality, path $\mathcal{P}_{uv}$ is expressed as
$$\mathcal{P}_{uv}:=uu^{1}_{i}s^{1}_{j}s\dots ww^{1}_{l}v^{1}_{n}v.$$
This leads to an equality
\begin{equation}\label{eqa:MF-3-1-4-4}
d^{\dagger}_{uv^{1}_{n}}+d^{\dagger}_{uw^{1}_{l}}+d^{\dagger}_{vu^{1}_{i}}+d^{\dagger}_{vs^{1}_{j}}=4d^{\dagger}_{uv}-6.
\end{equation}
Armed with the results in Eqs.(\ref{eqa:MF-3-1-4-3}) and (\ref{eqa:MF-3-1-4-4}), quantity $\mathcal{W}_{t}(4)$ is rewritten as
\begin{equation}\label{eqa:MF-3-1-4-5}
\begin{aligned}\mathcal{W}_{t}(4)&:=4[\mathcal{W}_{t}(1)-3(|\mathcal{T}_{t-1}|-1)]-6\left(
                                                             \begin{array}{c}
                                                               |\mathcal{T}_{t-1}|-1 \\
                                                               2 \\
                                                             \end{array}
                                                           \right)+4(|\mathcal{T}_{t-1}|-1)\\
&=4\mathcal{W}_{t}(1)-(|\mathcal{T}_{t-1}|-1)(3|\mathcal{T}_{t-1}|+2).
\end{aligned}
\end{equation}
Taking into account a fact that both terms in the second line of Eq.(\ref{eqa:MF-3-1-4-5}) are easy to answer, the closed-form solution of $\mathcal{W}_{t}(4)$ is more conveniently calculated in view of Eq.(\ref{eqa:MF-3-1-4-5}). It is clear to the eye that such a transformation from Eq.(\ref{eqa:MF-3-1-4-2}) to Eq.(\ref{eqa:MF-3-1-4-5}) is quite effective. At the same time, it should be mentioned that the light shed by the development of Eq.(\ref{eqa:MF-3-1-4-5}) will be proved helpful in the following derivations.

\emph{Stage 5.} Let us divert our insight into discussion about distance $d^{\dagger}_{u^{1}_{i}v^{1}_{j}}$ between two vertices, say $u^{1}_{i}$ and $v^{1}_{j}$, from sets $A_{u}(t)$ and $A_{v}(t)$. It is worth noticing that subscript $u$ can be the same as $v$. In other words, the result in Eq.(\ref{subeq:MF-3-1-3-3}) is covered by the following formula. By definition, given a couple of distinct sets $A_{u}(t)$ and $A_{v}(t)$, the summation $\mathcal{W}_{t}(A_{u}(t);A_{v}(t))$ that is defined as distances $d^{\dagger}_{u^{1}_{i}v^{1}_{j}}$ of all possible vertex pairs is given by

\begin{equation}\label{eqa:MF-3-1-5-1}
\begin{aligned}\mathcal{W}_{t}(A_{u}(t);A_{v}(t))&:=k_{u}(t-1)k_{v}(t-1)d^{\dagger}_{uv}-2+2(k_{u}(t-1)-1)(k_{v}(t-1)-1)
\end{aligned}
\end{equation}
where $A_{u}(t)$ is distinct with $A_{v}(t)$ as required.

To make further progress, from Eqs.(\ref{subeq:MF-3-1-3-3}) and (\ref{eqa:MF-3-1-5-1}), we can write

\begin{equation}\label{eqa:MF-3-1-5-2}
\begin{aligned}\mathcal{W}_{t}(5)&:=\frac{1}{2}\sum_{u,v\in\mathcal{T}_{t-1}}\mathcal{W}_{t}(A_{u}(t);A_{v}(t))+\sum_{v\in\mathcal{T}_{t-1}}\mathcal{W}_{t}(A_{v}(t))\\
&=\frac{1}{2}\sum_{u,v\in\mathcal{T}_{t-1}}k_{u}(t-1)k_{v}(t-1)d^{\dagger}_{uv}-\left(
                                                             \begin{array}{c}
                                                               |\mathcal{T}_{t-1}| \\
                                                               2 \\
                                                             \end{array}
                                                           \right)\\
&\quad+\sum_{u,v\in\mathcal{T}_{t-1}}(k_{u}(t-1)-1)(k_{v}(t-1)-1)\\
&\quad+\sum_{v\in\mathcal{T}_{t-1}}[k_{v}(t-1)]^{2}-2(|\mathcal{T}_{t-1}|-1).
\end{aligned}
\end{equation}

In general, determining solution to quantity $\mathcal{W}_{t}(5)$ according to Eq.(\ref{eqa:MF-3-1-5-2}) is a tough task. To address this issue, we still attempt to take advantage of transformational method. To be more concrete, the detailed calculations are deduced based on a similar manner to that used in stage 4. That is to say, we continue to focus on path $\mathcal{P}_{uv}$ and also make use of some notations mentioned above. Precisely, there are two sub-cases contained in the following discussions. For path $\mathcal{P}_{uv}$ with length $3$, distance $d^{\dagger}_{u^{1}_{i}v^{1}_{j}}$ is calculated to equal $1$. Additionally, if the length of path $\mathcal{P}_{uv}$ is greater than $3$, we have an equality like that in Eq.(\ref{eqa:MF-3-1-4-4}), which is given by

\begin{equation}\label{eqa:MF-3-1-5-3}
d^{\dagger}_{u^{1}_{i}v^{1}_{n}}+d^{\dagger}_{u^{1}_{i}w^{1}_{l}}+d^{\dagger}_{s^{1}_{j}v^{1}_{n}}+d^{\dagger}_{s^{1}_{j}w^{1}_{l}}=4d^{\dagger}_{uv}-12.
\end{equation}
Analogously, we can re-express quantity $W_{t}(5)$ in terms of the next equation

\begin{equation}\label{eqa:MF-3-1-5-4}
\begin{aligned}\mathcal{W}_{t}(5)&:=4[\mathcal{W}_{t}(1)-3(|\mathcal{T}_{t-1}|-1)]-12\left(
                                                             \begin{array}{c}
                                                               |\mathcal{T}_{t-1}|-1 \\
                                                               2 \\
                                                             \end{array}
                                                           \right)+(|\mathcal{T}_{t-1}|-1)\\
&=4\mathcal{W}_{t}(1)-(|\mathcal{T}_{t-1}|-1)(6|\mathcal{T}_{t-1}|-1).
\end{aligned}
\end{equation}

\emph{Stage 6.} Here, we pay more attention on determination of distance $d^{\dagger}_{uv^{2}_{i}}$. Note that vertex $u$ is different from $v$. Along the same line as in stage 4, the summation $\mathcal{W}_{t}(u;B_{v}(t))$ over distances between vertex $u$ and each vertex $v^{2}_{i}$ in set $B_{v}(t)$ is calculated to yield

\begin{equation}\label{eqa:MF-3-1-6-1}
\begin{aligned}\mathcal{W}_{t}(u;B_{v}(t))&:=\sum_{v(\neq u)\in\mathcal{T}_{t-1}}\sum_{v^{2}_{i}\in B_{v}(t)}d^{\dagger}_{uv^{2}_{i}}\\
&=\sum_{v(\neq u)\in\mathcal{T}_{t-1}}\left[k_{v}(t-1)d^{\dagger}_{uv}+k_{v}(t-1)\right].
\end{aligned}
\end{equation}
In order to simplify the next derivation, we now build up a bridge between quantities $\mathcal{W}_{t}(u;B_{v}(t))$ and $\mathcal{W}_{t}(u;A_{v}(t))$, which is
\begin{equation}\label{eqa:MF-3-1-6-2}
\mathcal{W}_{t}(u;B_{v}(t)):=\mathcal{W}_{t}(u;A_{v}(t))+\sum_{v(\neq u)\in\mathcal{T}_{t-1}}2.
\end{equation}

And then, we have

\begin{equation}\label{eqa:MF-3-1-6-3}
\begin{aligned}\mathcal{W}_{t}(6)&:=\sum_{u\in\mathcal{T}_{t-1}}\mathcal{W}_{t}(u;B_{v}(t))\\
&=\sum_{u\in\mathcal{T}_{t-1}}\sum_{v(\neq u)\in\mathcal{T}_{t-1}}k_{v}(t-1)d^{\dagger}_{uv}+\sum_{u\in\mathcal{T}_{t-1}}\sum_{v(\neq u)\in\mathcal{T}_{t-1}}k_{v}(t-1)\\
&=\mathcal{W}_{t}(4)+\sum_{u\in\mathcal{T}_{t-1}}\sum_{v(\neq u)\in\mathcal{T}_{t-1}}2\\
&=\mathcal{W}_{t}(4)+2|\mathcal{T}_{t-1}|(|\mathcal{T}_{t-1}|-1).
\end{aligned}
\end{equation}

With Eq.(\ref{eqa:MF-3-1-4-5}), we can obtain

\begin{equation}\label{eqa:MF-3-1-6-4}
\mathcal{W}_{t}(6):=4\mathcal{W}_{t}(1)-(|\mathcal{T}_{t-1}|-1)(|\mathcal{T}_{t-1}|+2).
\end{equation}

\emph{Stage 7.} For distance $d^{\dagger}_{u^{2}_{i}v^{2}_{j}}$ between two vertices, say $u^{2}_{i}$ and $v^{2}_{j}$, from distinct sets $B_{u}(t)$ and $B_{v}(t)$, using a similar deduction to that used to develop quantity $\mathcal{W}_{t}(A_{u}(t);A_{v}(t))$ in stage 5 produces

\begin{equation}\label{eqa:MF-3-1-7-1}
\begin{aligned}\mathcal{W}_{t}(B_{u}(t);B_{v}(t))&:=k_{u}(t-1)k_{v}(t-1)d^{\dagger}_{uv}+2k_{u}(t-1)k_{v}(t-1).
\end{aligned}
\end{equation}

As above, we also find a close connection of quantity $\mathcal{W}_{t}(B_{u}(t);B_{v}(t))$ to $\mathcal{W}_{t}(A_{u}(t);A_{v}(t))$, namely,

\begin{equation}\label{eqa:MF-3-1-7-2}
\begin{aligned}\mathcal{W}_{t}(B_{u}(t);B_{v}(t))&:=\mathcal{W}_{t}(A_{u}(t);A_{v}(t))+2(k_{u}(t-1)+k_{v}(t-1)).
\end{aligned}
\end{equation}
After that, with the help of Eqs.(\ref{eqa:MF-3-1-5-4}) and (\ref{eqa:MF-3-1-7-2}), we obtain

\begin{equation}\label{eqa:MF-3-1-7-3}
\begin{aligned}\mathcal{W}_{t}(7)&:=\frac{1}{2}\sum_{u,v\in\mathcal{T}_{t-1}}\mathcal{W}_{t}(B_{u}(t);B_{v}(t))+\sum_{v\in\mathcal{T}_{t-1}}\mathcal{W}_{t}(B_{v}(t))\\
&=\mathcal{W}_{t}(5)+\sum_{u\in\mathcal{T}_{t-1}}\sum_{v(\neq u)\in\mathcal{T}_{t-1}}(k_{u}(t-1)+k_{v}(t-1))\\
&=\mathcal{W}_{t}(5)+\sum_{u\in\mathcal{T}_{t-1}}\left[(|\mathcal{T}_{t-1}|-2)k_{u}(t-1)+2(|\mathcal{T}_{t-1}|-1)\right]\\
&=4\mathcal{W}_{t}(1)-(|\mathcal{T}_{t-1}|-1)(2|\mathcal{T}_{t-1}|+3).
\end{aligned}
\end{equation}
Notice also that Eq.(\ref{subeq:MF-3-1-3-4}) has been used.

\emph{Stage 8.}  The last part of our tasks is to measure distance $d^{\dagger}_{u^{1}_{i}v^{2}_{j}}$ of vertices $u^{1}_{i}$ and $v^{2}_{j}$. This may be finished in the following two cases. When both vertices are connected to an identical vertex $v$, it is clear to find an equality

\begin{equation}\label{eqa:MF-3-1-8-1}
\begin{aligned}\mathcal{W}_{t}(A_{v}(t);B_{v}(t)&:=\sum_{v^{1}_{i}\in A_{v}(t)}\sum_{v^{2}_{i}\in B_{v}(t)}d^{\dagger}_{v^{1}_{i}v^{2}_{j}}\\
&=2[k_{v}(t-1)]^{2}.
\end{aligned}
\end{equation}
And then, we establish

\begin{equation}\label{eqa:MF-3-1-8-2}
\begin{aligned}\mathcal{W}_{t}(8)&:=\sum_{v\in\mathcal{T}_{t-1}}\mathcal{W}_{t}(A_{v}(t);B_{v}(t))\\
&=2\sum_{v\in\mathcal{T}_{t-1}}[k_{v}(t-1)]^{2}.
\end{aligned}
\end{equation}

On the other hand, if vertices $u^{1}_{i}$ and $v^{2}_{j}$ are from distinct sets $A_{u}(t)$ and $B_{v}(t)$, respectively, then we build up

\begin{equation}\label{eqa:MF-3-1-8-3}
\begin{aligned}\mathcal{W}_{t}(A_{u}(t);B_{v}(t))&:=\sum_{u^{1}_{i}\in A_{u}(t)}\sum_{v^{2}_{i}\in B_{v}(t)}d^{\dagger}_{u^{1}_{i}v^{2}_{j}}\\
&=k_{u}(t-1)k_{v}(t-1)d^{\dagger}_{uv}+2k_{u}(t-1)(k_{v}(t-1)-1)\\
&=\mathcal{W}_{t}(B_{u}(t);B_{v}(t))-2k_{v}(t-1),
\end{aligned}
\end{equation}
where we have used Eq.(\ref{eqa:MF-3-1-7-1}). Also, based on Eq.(\ref{eqa:MF-3-1-7-3}), it is easy to check the following equation

\begin{equation}\label{eqa:MF-3-1-8-4}
\begin{aligned}\mathcal{W}_{t}(9)&:=\sum_{u,v\in\mathcal{T}_{t-1}}\mathcal{W}_{t}(A_{u}(t);B_{v}(t))+\sum_{v\in\mathcal{T}_{t-1}}\mathcal{W}_{t}(A_{v}(t);B_{v}(t))\\
&=2\left(\mathcal{W}_{t}(7)-\sum_{v\in\mathcal{T}_{t-1}}\mathcal{W}_{t}(B_{v}(t))\right)-2\sum_{u\in\mathcal{T}_{t-1}}\sum_{v(\neq u)\in\mathcal{T}_{t-1}}k_{v}(t-1)+2\sum_{v\in\mathcal{T}_{t-1}}[k_{v}(t-1)]^{2}\\
&=2\mathcal{W}_{t}(7)-2\left(\sum_{v\in\mathcal{T}_{t-1}}[k_{v}(t-1)]^{2}-2(|\mathcal{T}_{t-1}|-1)\right)-4(|\mathcal{T}_{t-1}|-1)^{2}+2\sum_{v\in\mathcal{T}_{t-1}}[k_{v}(t-1)]^{2}\\
&=2\mathcal{W}_{t}(7)-4(|\mathcal{T}_{t-1}|-1)^{2}+4(|\mathcal{T}_{t-1}|-1)\\
&=8\mathcal{W}_{t}(1)-2(|\mathcal{T}_{t-1}|-1)(4|\mathcal{T}_{t-1}|-1).
\end{aligned}
\end{equation}

Until now, we have exhaustively enumerated all possible combinatorial ways when considering distance of an arbitrary pair of vertices in tree $\mathcal{T}_{t}$. After that, putting all things together, we can see that, after performing some simple arithmetics, Wiener index $\mathcal{W}_{t}$ of scale-free tree network $\mathcal{T}_{t}$ satisfies
\begin{equation}\label{eqa:MF-3-1-8-5}
\mathcal{W}_{t}=\sum_{i=1}^{9}\mathcal{W}_{t}(i)-\mathcal{W}_{t}(3)-\mathcal{W}_{t}(8)=75\mathcal{W}_{t-1}-20|\mathcal{T}_{t-1}|^{2}+20|\mathcal{T}_{t-1}|.
\end{equation}
This means that Wiener index $\mathcal{W}_{t}$ of tree network $\mathcal{T}_{t}$ is precisely estimated in an iterative manner. Then, using the elegant relation in Eq.(\ref{eqa:MF-3-0-4}), the closed-from solution of mean hitting time $\langle\mathcal{H}_{t}\rangle$ is derived. This is complete. \qed

\end{document}